\documentclass[12pt]{article}

\usepackage{amsfonts}
\usepackage{amssymb}
\usepackage{amscd}
\usepackage{verbatim}
\usepackage{latexsym}
\font\rurm=wncyr10 scaled \magstep1

\begin{document}

\title{Constructing semisimple $p$-adic Galois representations with prescribed properties}

\author{Chandrashekhar Khare, University of Utah and 
TIFR
\thanks{some of the work on this paper was done during a visit to Universit\'e Paris 7 which was supported 
by Centre franco-indien pour la promotion de la recherche avanc{\'e}e (CEFIPRA) under
Project 2501-1 {\it Algebraic Groups in Arithmetic and Geometry}}, \\
Michael Larsen, Indiana University 
\thanks{partially supported
by NSF Grant DMS-0100537}, \\
\& Ravi Ramakrishna, Cornell University and McGill University 
\thanks{partially supported
by NSF Grant DMS-0102173 and the AMS 
Centennial Research Fellowship. The author thanks UC-Berkeley
and Princeton University for their hospitality during visits in 2003.}}

\date{}

\maketitle
\newtheorem{theorem}{Theorem}
\newtheorem{lemma}[theorem]{Lemma}
\newtheorem{prop}[theorem]{Proposition}
\newtheorem{fact}[theorem]{Fact}
\newtheorem{cor}[theorem]{Corollary}
\newtheorem{example}{Example}
\newtheorem{conj}[theorem]{Conjecture}
\newtheorem{definition}[theorem]{Definition}
\newtheorem{quest}[theorem]{Question}
\newtheorem{ack}{Acknowledgemets}
\newcommand{\rhobar}{\overline{\rho}}
\newcommand{\Ad}{{\rm Ad}}

\def\sha{{{\textnormal{\rurm{Sh}}}}}
\def\qrhob{\bf Q(\bar{\rho})}
\def\eps{\epsilon}
\def\rhobar{ {\bar {\rho} } }
\def\wfq{W({\bf{F_q}})}
\def\ad{Ad^0\bar{\rho}}
\def\adst{ Ad^0\bar{\rho}(1)}
\def\gal{Gal( \bar{ {\mathbb Q}}/{\mathbb Q})}
\newcommand{\rhob}{\overline{\rho}}

\section{Introduction and sketch of proof of the main theorem}

The study of $p$-adic representations of absolute Galois groups of 
number fields, i.e., continuous representations 
${\rho}:G_{\bf K} \rightarrow GL_n({\mathbb Q}_p)$ 
with $G_{\bf K}$ the absolute Galois group of a number field and $p$ a prime, is one of the central themes of modern number theory.
The ones studied the most are those which arise from the \'etale cohomology of smooth, projective varieties
over number fields. These have the 
striking property (due to Weil, Dwork, Grothendieck et al.)
that they are ramified at finitely many primes and are {\it rational} over a number field ${\bf L}$, i.e., for all
primes $r$ that are not ramified in $\rho$, the characteristic polynomial attached to the conjugacy class
of ${\rm Frob}_r$ in the image of $\rho$ (with ${\rm Frob}_r$ the Frobenius substitution at $r$) has coefficients in ${\bf L}[X]$.
Each belongs to a {\it compatible} 
family of Galois representations and is {\it pure} of some weight $k$.
In this paper we give a purely Galois theoretic method for constructing
semisimple continuous representations ${\rho}:G_{\bf K} \rightarrow GL_n({\mathbb Q}_p)$
that are {\it density $1$ rational} over ${\mathbb Q}$ in the 
critical case of ${\bf K}={\mathbb Q}$ and $n=2$ and {\it density $1$ pure} 
(see Definition
\ref{purity} below). Unfortunately these representations are ramified
at infinitely many primes: even then being rational
is a strong condition as the primes that 
are ramified in a semisimple representation $\rho$ 
have density 0 (see [Kh-Raj]
and also Definition \ref{purity} below).
As typically our constructions give 
infinitely ramified representations, we get examples
of density $1$ 
rational $p$-adic representations that are {\it emphatically} non-geometric although they do arise often as $p$-adic limits
of geometric representations. In fact  for $p \geq 5$
any representation
surjective onto $GL_2({\mathbb Z}/p{\mathbb Z})$ lifts to 
${\rho}:G_{\mathbb Q} \rightarrow GL_2({\mathbb Q}_p)$
which is rational over ${\mathbb Q}$.
Our methods also 
allow us to also construct density $1$ compatible lifts of 
almost any given pair of 2-dimensional mod $p$ and mod $q$ representations 
of $G_{\mathbb Q}$ (see Definition \ref{compatible} below). 
We describe below the main new method of this paper.

We start with a residual representation
$\rhob:\gal \rightarrow GL_2({\mathbb Z}/p{\mathbb Z})$ and assume $p \geq 5$.
Our aim is to lift $\rhob$ to $GL_2({\mathbb Z}_p)$. In many cases, this
was done in [R3].
The strategy there (and here) was to successively lift from
$GL_2({\mathbb Z}/p^{m}{\mathbb Z})$ to $GL_2({\mathbb Z}/p^{m+1} {\mathbb Z})$.
Auxiliary primes at which ramification was allowed were introduced. First it was
shown that any global obstructions to lifting
would be realised locally.
Then it was shown that if there was a local obstruction to lifting,
one could choose a different lift to 
$GL_2({\mathbb Z}/p^m{\mathbb Z})$, congruent mod $p^{m-1}$ to the old one
(using elements of  a global $H^1$),
where all local obstructions vanished.
(Though see [T2] where these two conditions are
handled simultaneously.)
The limiting  representation
to $GL_2({\mathbb Z}_p)$ was ramified at a finite set of primes.
The main difficulty in [R3] was as follows:
we introduced more primes
of ramification to make the global $H^1$ bigger, 
with the intent of allowing us to use more
global $H^1$ elements so that we could remove
local obstructions to lifting. As we did this we introduced more
potential local obstructions at the new ramified 
primes. In most cases these competing
forces were exactly balanced against each other and a unique lift (for the given
ramified set) was shown to exist. In the cases that are not
dealt with in [R3] the methods there were not able to balance the
competing forces.

Here our approach is somewhat different. In addition to producing a lift
at each stage, we want to choose, at our discretion, specified characteristic
polynomials of Frobenius at each stage
of any finite number of unramified primes. 
(There will be more and more characteristic polynomials that we choose
at each stage.)
Thus
at each stage of the deformation process
we have {\em far more} local conditions to arrange than global
$H^1$ elements to adjust them by with ramification allowed only at the auxiliary set that occurred at the earlier stage. 
We are able to do all this 
only by allowing our $GL_2({\mathbb Z}_p)$ representation
to be ramified at infinitely many primes.
Of course its mod $p^m$ reduction, for any $m$ is ramified at only finitely
many primes. At each stage of the lifting process we will impose more
local conditions, and require ramification at more primes.

\vskip1em\noindent
{\bf Main Theorem } {\em Let $p \geq 5$ and
$\rhob:\gal \rightarrow GL_2({\mathbb Z}/p{\mathbb Z})$ be given such that
the image of $\rhob$ contains $SL_2({\mathbb Z}/p{\mathbb Z})$.
If $p=5$  assume further $\rhob$ is surjective.
Assume det$\rhob=\chi^k$ where $\chi$
is the mod $p$ cyclotomic character and $1 \leq k \leq p-1$.
Then there exists a deformation $\rho$ of $\rhob$ to ${\mathbb Z}_p$
such that 
$\rho|_{G_p}$ is potentially semistable,
$\rho$ is unramified at a density $1$ set of primes $R$,
and for all but finitely many unramified primes $r$
the characteristic polynomial of Frobenius at $r$
is in ${\mathbb Z}[x]$  pure of weight $k$.} 
%(in the sense of Definition \ref{purity} below).
%Furthermore, for the density $0$ infinite set of primes which we added
%to our ramification set, we may choose the local representation
%at any of these primes to be whatever we want, provided our choice
%is consistent with the local representation at the mod $p^m$ level
%these primes were introduced. In particular, if $\rhob|_{G_p}$
%is unramified, we may choose $\rho|_{G_p}$ to be unramified.

\vskip1em\noindent
Remarks: 
1) In the main theorem above we can also ensure that
the lift $\rho$ is ramified at  infinitely many primes. 
Getting a pure, rational lift such as 
$\rho$ above that is finitely ramified by using only Galois
theoretic methods seems extremely hard, if not impossible.
(Though see [T1] for a result in this direction using 
geometric techniques.)
In fact given $\rhob$ as above, 
and fixing the characteristic polynomial $f_r(X)$ at even one
prime $r$ that is consistent with $\rhob$, 
we see no way of getting a potentially semistable at $p$ lift $\rho$ such 
that $\rho$ is unramified at $r$ and has characteristic polynomial 
$f_r(X)$ with $\rho$ finitely ramified!
\newline\noindent
2) The method used to prove the 
theorem above is quite involved and we believe to be of
independent interest. It forms the core of 
the paper (see the sketch below and subsection 2.1).

\vskip1em\noindent

We point out in passing 2 amusing consequences of the method of the proof of the main theorem:
\begin{itemize}
 \item We can construct a continuous, semisimple representation $\rho:G_{\mathbb Q} \rightarrow GL_2({\mathbb C}_p)$
 with ${\mathbb C}_p$ the completion of $\overline{\mathbb Q}_p$ 
that is not conjugate to a representation
into $GL_2(\overline{\mathbb Q}_p)$ 
(we can even get such representations that are unramified at 
a density 1 set of primes
and at these primes the characteristic polynomials of the 
Frobenii are defined over $\overline{\mathbb Q}$ regarded inside
${\mathbb C}_p$)! The method also 
gives a way to lift surjective mod $p$ representations
$\rhobar:G_{\mathbb Q} \rightarrow GL_2({\mathbb Z}/p{\mathbb Z})$ 
to representations into $GL_2({\bf K})$ with 
${\bf K}$ {\it any} fixed finite extension
of ${\mathbb Q}_p$ with field of definition of the 
lift being ${\bf K}$. 

  \item We can construct a {\it surjective}, 
continuous representation 
$\rho:G_{\mathbb Q} \rightarrow GL_2({\bf K})$ with ${\bf K}$ 
a finite extension
of ${\mathbb Q}_p$ such that there is a finite 
extension ${\bf F}/{\mathbb Q}$ (regarded as embedded in ${\bf K}$) such 
that the characteristic
polynomials at all unramified primes split over ${\bf F}$ (these 
are called ${\bf F}$-split representations in [Kh3]: the existence of such a 
$\rho$ also shows that 
for the question at end of [Kh3] to have a positive 
answer it is necessary to restrict attention to finitely 
ramified representations).  It is not hard to 
show (for instance using the arguments of [Kh3]) 
that for such a $\rho$, for any prime $q$
that splits in ${\bf F}$, there cannot be a 
representation $\rho':G_{\mathbb Q} \rightarrow 
GL_2(\overline{\mathbb Q}_q)$ that is
density one compatible with $\rho$.
\end{itemize}

\noindent
{\bf Sketch of proof of Main Theorem:}
Starting with $\rhob$, we let $S$ be a finite
set of primes containing $p$, those ramified in $\rhob$, and 
enough of what we call {\em nice} primes for $\rhob$
(see Definition~\ref{nicedef}) so that
global obstructions to lifting to $GL_2({\mathbb Z}/p^2{\mathbb Z})$
can be detected locally. 
Fix the determinant of all our lifts to be $\epsilon^k$
where $k$ is chosen suitably large.

Let $S_2=S$.
Once and for all, for each $v \in S_2$ choose a potentially semistable
local deformation
of $\rhob|_{G_v}$ to $GL_2({\mathbb Z}_p)$.
(This is not difficult. See [R3] for instance.
The condition of potential semistability is vacuous for $v \neq p$.)
Deform $\rhob$ to $\rho_2:G_S
\rightarrow GL_2({\mathbb Z}/p^2{\mathbb Z})$. 
Note we do not know $\rho_2|_{G_v}$ is the mod $p^2$ reduction
of our preselected deformation to ${\mathbb Z}_p$.
Now we let $R_2$ consist of all primes
beneath a certain (large) bound that are not in $S_2$. 
Once and for all choose characteristic polynomials in ${\mathbb Z}[x]$
for
all $r \in R_2$.

Using Lemma~\ref{unram}
we find a collection $Q_2$ of $\rho_2$-nice primes
such that the map 
\begin{equation}H^1(G_{S_2 \cup Q_2},\ad) \rightarrow
\oplus_{v \in S_2} H^1(G_v,\ad) \oplus_{v \in R_2}H^1_{nr}(G_r,\ad)
\end{equation}
is an isomorphism. This implies there is a
unique $f_2 \in H^1(G_{S_2 \cup Q_2},\ad)$
such that $(I+pf_2)\rho_2 $ is the mod $p^2$ reduction
of our preselected deformation to ${\mathbb Z}_p$ for all $v \in S_2$
{\em and} the characteristic polynomials at all primes in $R_2$
are the mod $p^2$ reductions of the preselected 
characteristic polynomials
as well. We may, however, have introduced local obstructions
to lifting at primes in
$Q_2$. We remove these obstructions (see Proposition \ref{polarisation}) by allowing more ramification
at a set $V_2$ of $\rho_2$-nice primes, that is, by adjusting by an element
of $H^1(G_{S_2 \cup Q_2 \cup V_2},\ad)$.
This set $V_2$ will have cardinality up to 
twice that of $Q_2$, but in the end there  will be
no obstructions to lifting at primes of 
$S_2 \cup Q_2 \cup V_2$, nor will we change anything
at primes in $S_2 \cup R_2$. {\em The existence of the set $V_2$ (and later
the sets $V_n$) is the key technical innovation required to prove the
Main Theorem.}

Put $S_3 = S_2 \cup Q_2 \cup V_2$.
Then we deform this last representation to 
$\rho_3:G_{S_3} \rightarrow GL_2({\mathbb Z}/p^3{\mathbb Z})$.
Once and for all  preselect lifts of $\rho_3|_{G_v}$
for all $v \in S_3 \backslash S_2$ to $GL_2({\mathbb Z}_p)$. 
%(Since such $v$ are $\rho_2$-nice, if $\rho_3$ is unramified at
%$v$ we may preselect a ramified or unramified lift as we please.)
Take $R_3$ to be the union of $R_2$ and all primes below
a certain bound not in $S_3$. For the primes in $R_3 \backslash R_2$,
once and for all choose pure weight $k$ characteristic polynomials
in ${\mathbb Z}[x]$ consistent with the mod $p^2$ reduction
of $\rho_3$. That we can guarantee these characteristic polynomials
have both roots with absolute value $p^{k/2}$ will follow from the fact
that the primes of $R_3 \backslash R_2$ are suitably big.
Now simply repeat the process.

In the limit we have a deformation of $\rhob$ to $GL_2({\mathbb Z}_p)$
that is, by [Kh-Raj], ramified at a density $0$ set of primes, 
and whose characteristic polynomials of Frobenius at all but
finitely many of the density $1$ set of
unramified primes are pure of weight $k$.

Serre has asked if the method of proof of 
the Main Theorem can be used to give another approach to Shafarevich's
result that realises any finite solvable group as the Galois group of
a finite extension of a given number field. We also remark that our techniques 
have some resemblance to a key ingredient in the proof of Shafarevich's theorem
in [NSW] (see Theorem 9.5.9 of loc. cit.).

We conclude the paper by
proving a result
about representations ramified at infinitely many primes
that shows that the main theorem
of [Kh-Raj] is essentially best possible (answering a
question of Serre).

\section{The main results}

\subsection{The toolbox}
\noindent
We make some preliminary observations.
\newline\noindent
1) In this paper we always deal with the cohomology of
$\ad$, the set of $2 \times 2$ trace zero
matrices over ${\mathbb Z}/p{\mathbb Z}$,
as opposed to the cohomology of $Ad\rhob$ the group
of all $2 \times 2 $ matrices. This basically means we are fixing
all determinants of all of our global
deformations once and for all.
We assume  det$\rhob=\chi^k$  where $\chi$ is the mod $p$
cyclotomic character and $0 \leq k \leq p-1$. We will always choose the
determinant of all of our deformations to be $\epsilon^k$
where $\epsilon$ is the $p$-adic cyclotomic
character. (Note $k$ is only well defined mod $p-1$, but we fix it once
and for all at the beginning of our construction.)
\newline\noindent
2) Let $\rho_{m}:G \rightarrow GL_n({\mathbb Z}/p^{m}{\mathbb Z})$
and  let $\rho_{m+1}:G \rightarrow GL_n({\mathbb Z}/p^{m+1} {\mathbb Z})$
be lifts 
of $\rho_{m}$ with both lifts of $\rhob$. Let $f \in H^1(G,\ad)$. Then
using the cocycle relation, one easily sees that the map
\begin{equation}g \mapsto (I+p^{m}f(g))\rho_{m+1}(g)\end{equation} 
is another lift of $\rho_m$ to $GL_2({\mathbb Z}/p^{m+1})$.
For $f$ a coboundary one gets that $ (I+p^{m}f(g))\rho_{m+1}(g)$
is the same {\em deformation} as $\rho_m$. See [M] for
this definition and an introduction to deformation theory.
If the centraliser of the image of $\rhob$ is exactly the scalars
then the deformations of $\rho_m$ to ${\mathbb Z}/p^{m+1} {\mathbb Z}$
form a principal homogeneous space over $H^1(G,\ad)$.
\newline\noindent
3) It is known that {\em all} $\rhob|_{G_v}$ admit deformations
to ${\mathbb Z}_p$. For $v \neq p$ this is mentioned in [R3], following
work of Diamond and Taylor.
If $v=p$ and $\rhob|_{G_p}$ is ramified, there is always a potentially
semistable deformation to ${\mathbb Z}_p$. This is worked out in [R3].
See also [B\"{o}].
For example, for $v=p$ and $\rhob|_{G_p}$ unramified
it is an easy exercise to show there are unramified deformations
to ${\mathbb Z}_p$. If $( \# \rhob(G_p),p)=1$ there
is the  Teichm\"{u}ller deformation of $\rhob|_{G_p}$ whose image is
isomorphic to $\rhob(G_p)$ and therefore 
of finite order and thus clearly potentially semistable.
\newline\noindent
4) As we lift from mod $p^m$ to mod $p^{m+1}$ we
will be choosing more primes whose characteristic polynomials we
choose once and for all at each stage.
If at the mod $p^m$ stage, we have a prime $r$ whose
characteristic polynomial we need to choose once and for all (consistent
with the mod $p^m$ deformation of course), then different choices
will lead to different 
mod $p^{m+1} {\mathbb Z}$ deformations, and we will be
restricted in our choice of characteristic polynomials for future
primes. Thus, in the end, when we have our deformation to ${\mathbb Z}_p$
we cannot change the characteristic polynomial(s) of any one (or finite number)
of primes. It is important to note we are {\em not} choosing all our
characteristic polynomials of unramified primes at the beginning.
\newline\noindent
5) In this paper, we are providing examples of certain pathologies.
We have not sought maximal generality. In particular we work only with
residual representations to $GL_2({\mathbb Z}/p{\mathbb Z})$ and deformations
to $GL_2({\mathbb Z}/p^m{\mathbb Z})$, 
rather than with residual representations
to $GL_n({\bf k})$ where ${\bf k}$ is a finite field.
We also assume det$\rhob$ is a power of the mod $p$ cyclotomic character,
the image of $\rhob$ contains $SL_2({\mathbb Z}/p{\mathbb Z})$, and that
$H^1(Im\,\rhob,Ad)=0$.
We expect it is
not difficult to extend	 our results to the general case.
\newline\noindent
6) In [R3] the parity of $\rhob$ played an important role. For $\rhob$ odd,
one had two `extra degrees of freedom', which were just enough to 
(usually) provide
a characteristic zero deformation that was potentially semistable. Here,
since we are circumventing the problem of having
more local conditions than global degrees of freedom, the parity of
$\rhob$ plays no role.
\begin{definition}\label{nicedef} 
Suppose $\rhob:\gal \rightarrow GL_2({\mathbb Z}/p{\mathbb Z})$ is given as in
the Main Theorem.
We say a prime $q$ is {\em nice} (for $\rhob$)
if 
\begin{itemize}
  \item $q$ is {\em not} $\pm 1$ mod $p$, 
  \item $\rhob$ is unramified at $q$,
  \item  the eigenvalues of $\rhob(\sigma_q)$ (where $\sigma_q$ is Frobenius
         at $q$) have ratio $q$. 
\end{itemize}
Let $\rho_m$ be a
deformation of $\rhob$ to $GL_2({\mathbb Z}/p^m{\mathbb Z})$.
We say a  prime $q$ is $\rho_m$-{\em nice} if 
\begin{itemize}
  \item $q$ is nice for $\rhob$,
  \item $\rho_m$ is unramified at $q$, and the eigenvalues
        of $\rho_m(\sigma_q)$ have ratio $q$. 
Note that since $q$ is nice,
the mod $p^m$ characteristic polynomial of $\rho_m(\sigma_q)$ 
has distinct roots that are units, so the eigenvalues
of $\rho_m(\sigma_q)$ are well-defined in ${\mathbb Z}/p^m{\mathbb Z}$. 
\end{itemize}
\end{definition}
Remark: That $\rho_m$-nice primes exist follows
from Fact~\ref{disjointness}. 
\vskip1em
If $q$ is nice then any deformation of $\rhob|_{G_q}$ will be tamely
ramified. Since the Galois group over ${\mathbb Q}_q$ of the maximal
tamely ramified extension is generated by Frobenius $\sigma_q$ and
a  generator of tame inertia $\tau_q$ subject to
the relation $\sigma_q\tau_q{\sigma^{-1}_q}=\tau^q_q$, a
versal deformation is specified by the images of
$\sigma_q$ and $\tau_q$. We simply give them here.
See [R1] for more details.
The versal ring is ${\mathbb Z}_p[[A,B]]/(AB)$ and up to twist
\begin{equation}\sigma_q \mapsto \left(\begin{array}{cc} q(1+A) & 0\\0 & (1+A)^{-1}
\end{array}\right),\,\,\,\,
\tau_q \mapsto \left(\begin{array}{cc} 1 & B\\0 & 1
\end{array}\right).\end{equation}
We will be interested in 
deformations of $\rhob$ where $A \mapsto 0$.
An arbitrary mod $p^m$ local deformation of $\rhob|_{G_q}$ for $q$
nice will certainly be unobstructed if, up to twist
$\sigma_q \mapsto \left(\begin{array}{cc} q & 0\\0 & 1
\end{array}\right)$.
If on the other hand, up to twist,
 $\sigma_q \mapsto \left(\begin{array}{cc}
q(1+\alpha p^{m-1}) & 0\\0 & 1-\alpha p^{m-1}
\end{array}\right)$ where $\alpha \neq 0$
and $\tau_q$ has nontrivial image, then we
will {\em not} be able to deform this
mod $p^m$ representation to ${\mathbb Z}_p$. At all of our nice primes
our aim will be to adjust the (local) representation so 
$A \mapsto 0$.
Thus we will pay particular attention to the image of Frobenius
at nice primes.
\begin{definition}\label{unobstructed} For $\rho_m$ a deformation of $\rhob$ to
${\mathbb Z}/p^m{\mathbb Z}$ and $q$ a nice prime for $\rhob$ we call
$\rho_m|_{G_q}$ {\em unobstructed} for $q$ if, in the notation 
above, $A \mapsto 0$. 
\end{definition}
\noindent
Remark: There will be many nice primes which we will {\em not}
introduce into our ramification set. Since these primes
are unramified, their deformation problems have no
obstruction. In this paper we use the term {\em unobstructed}
almost exclusively for nice primes in the context of
Definition~\ref{unobstructed}.
\vskip1em

\begin{lemma}\label{nice}
For $M$ an unramified $G_q$-module, define $H^1_{nr}(G_q,M)$ to be the unramified
cohomology classes in $H^1(G_q,M)$.
Put $h^i=\dim(H^i(G_q,\ad))$ and $h^i_{d} =\dim(H^i(G_q,\adst))$
where $\adst:=Hom(\ad,\mu_p)$ is the ${\mathbb G}_m$-dual.
Then if $q$ is nice  we have
\begin{equation}h^0=1,h^1=2,h^2=1 \quad \,\,\,\,\,h^0_{d}=1, h^1_{d}=2,h^2_{d}=1\end{equation}
and 
\begin{equation}\dim(H^1_{nr}(G_q,\ad))
=1=\dim(H^1_{nr}(G_q,\adst)).\end{equation}
\end{lemma}
Proof: As a ${\mathbb Z}/p{\mathbb Z}[G_q]$ module
we have $\ad \simeq {\mathbb Z}/p{\mathbb Z} \oplus {\mathbb Z}/p{\mathbb Z}(1)
\oplus {\mathbb Z}/p{\mathbb Z}(-1)$ and
$\adst \simeq {\mathbb Z}/p{\mathbb Z}(1) \oplus {\mathbb Z}/p{\mathbb Z}
\oplus {\mathbb Z}/p{\mathbb Z}(2)$. 
From
local Galois cohomology, using that $p \geq 5$ and
$q$ is not $\pm 1$ mod $p$,
one easily sees that ${\mathbb Z}/p{\mathbb Z}(-1)$ and ${\mathbb Z}/p{\mathbb Z}(2)$
have trivial cohomology. So we only need to compute the cohomology
of $N:={\mathbb Z}/p{\mathbb Z} \oplus {\mathbb Z}/p{\mathbb Z}(1)$. 
This is a routine
computation using local duality and the local Euler characteristic.
We note 
$$H^1(G_q,\ad) \simeq H^1(G_q,\adst) \simeq H^1(G_q,N)$$
\begin{equation}
\simeq H^1(G_q,{\mathbb Z}/p{\mathbb Z}) 
\oplus H^1(G_q,{\mathbb Z}/p{\mathbb Z}(1)).\end{equation}
\hfill$\square$
\vskip1em\noindent
Remarks: 1) Suppose $\rho_m$ is a deformation of $\rhob$ to ${\mathbb Z}/p^m{\mathbb Z}$
and $q$ is a nice prime for $\rhob$. Let $f \in H^1(\gal,\ad)$.
When we speak of $f(\sigma_q)$ we will mean the
diagonal values of this trace zero matrix. 
Equivalently, since 
\begin{equation}H^1(G_q,\ad) \simeq H^1(G_q,{\mathbb Z}/p{\mathbb Z}) \oplus H^1(G_q,{\mathbb Z}/p{\mathbb Z}(1)),\end{equation}
$f(\sigma_q)$ is the value at $\sigma_q$ of the projection of $f$ on the first factor
in the direct summand. 
For $\phi \in H^1(\gal,\adst)$
we make a similar definition for $\phi(\sigma_q)$.
\newline\noindent
2) Let $\rho_m$ be a deformation of $\rhob$ as usual. Let $q$ be
nice for $\rhob$ and suppose
$\rhob_{m-1}|_{G_q}$ is a mod $p^{m-1}$ deformation
of $\rhob|_{G_q}$ which is unobstructed.
Suppose also that
in $\rho_m|_{G_q}$ that $A \not \mapsto 0$ and
$f \in H^1(\gal,\ad)$ is such that $f(\sigma_q) \neq 0$.
Then for some $\beta \in {\mathbb Z}/p{\mathbb Z}$ we have that
$(I+\beta f  p^{m-1})\rho_m|_{G_q}$ is unobstructed.
\vskip1em
We recall a proposition of Wiles (Prop. $1.6$ of [W])
which we adapt slightly for our purposes.

\begin{fact}\label{wiles1.6} Let $M$ be a finite dimensional
vector space over ${\mathbb Z}/p{\mathbb Z}$
with a $\gal$ action.
Let $S$ be a finite set of
primes containing $p$
and all primes that are ramified in the field fixed by the kernel
of the action of $\gal$ on $M$. For each $v \in S$ let 
${\cal L}_v \subset H^1(G_v,M)$ be a subspace with annihilator
${\cal L}^{\perp}_v \subset H^1(G_v,M(1))$ under the local pairing.
Define $H^1_{{\cal L}}(G_S,M)$ and $H^1_{{\cal L}^{\perp}}(G_S,M(1))$
to be, repectively, the kernels of the maps
\begin{equation}H^1(G_S,M) \rightarrow \oplus_{v \in S}
 \frac{H^1(G_v,M)}{{\cal L}_v},\,\,\,\,
H^1(G_S,M(1)) \rightarrow   \oplus_{v \in S}
\frac{H^1(G_v,M(1))}{{\cal L}^{\perp}_v}.\end{equation}
Then
$$\dim(H^1_{ {\cal L}}(G_S,M))
-\dim(H^1_{ {\cal L}^{\perp}}(G_S,M(1))) =$$
\begin{equation}\label{fact4} \dim(H^0(G_S,M))-\dim(H^0(G_S,M(1)))+
\sum_{v \in S} \left(\dim( {\cal L}_v )
- \dim(H^0(G_v,M)) \right).\end{equation}
\end{fact}
Proof: See Proposition $1.6$ of [W]  or $8.6.20$ of [NSW]. The result
follows from the long exact sequence of global Galois
cohomology.
\hfill $\square$
\vskip1em\noindent
Remark: The above groups are often called the Selmer and dual Selmer groups
for the set $S$ and local conditions ${\cal L}$ and ${\cal L}^{\perp}$
respectively.
In practice these groups are extremely difficult to compute as the class
groups of the fields fixed by the kernel of the Galois action on $M$ and $M(1)$
enter into the computations. However the formula shows the 
difference in dimension between the Selmer and dual Selmer groups
for a set of primes $S$ and local conditions ${\cal L}$ 
and ${\cal L}^{\perp}$ can be readily
computed. 
\vskip1em
We return to our set-up with $\rhob$.

\begin{fact}\label{poitoutate}
There is an exact sequence 
\begin{equation}H^1_{ {\cal L}_v}(G_S,\ad) \rightarrow
H^1(G_S,\ad) \rightarrow \oplus_{v \in S}\frac{H^1(G_v,\ad)}{{\cal L}_v}
 \rightarrow \end{equation}
 $$H^1_{ {\cal L}^{\perp}_v}(G_S,\adst)\rightarrow H^2(G_S,\ad) 
 \rightarrow \oplus_{v \in S} H^2(G_v,\ad).$$
\end{fact}
Proof: This follows from the Poitou-Tate exact sequence.
\hfill $\square$

\begin{fact}\label{disjointness}
Let $\rho_m$ be a deformation of $\rhob$ to ${\mathbb Z}/p^m{\mathbb Z}$
unramified outside $S$.
Let $\{f_1,...,f_n\}$  be
linearly independent in
$H^1(\gal,\ad)$
and $\{ \phi_1,...,\phi_r \}$ be 
linearly independent
in
$H^1(\gal,\adst)$.
Let $\mathbb{Q}(Ad^0(\rhobar))$ be the field
fixed by the kernel of the action of $G_{\mathbb{Q}}$
on $Ad^0(\rhobar)$. Let ${\bf K}=\mathbb{Q}(Ad^0(\rhobar),\mu_p)$
be the field obtained by adjoining the $p$th roots of unity to
$\mathbb{Q}(Ad^0(\rhobar))$.
We denote by ${\bf K}_{f_i}$ and ${\bf K}_{\phi_j}$
the fixed fields of the kernels of the restrictions
of $f_i$ and  $\phi_j$ to $G_{{\bf K}}$,
the absolute Galois group of ${\bf K}$.
Also, as ${\mathbb Z}/p{\mathbb Z}[Gal({\bf K}/{\mathbb Q})]$-modules,
$Gal({\bf K}_{f_i}/{\bf K})$ and $Gal({\bf K}_{\phi_j}/{\bf K})$
are isomorphic, respectively, to $\ad$ and $\adst$.
Let ${\bf P}_m$
be the fixed field of the kernel of the restriction of
the projectivisation of $\rho_m$ to $G_{{\bf K}}$.
Then each of the fields ${\bf K}_{f_i}$, ${\bf K}_{\phi_j}$
${\bf P}_m$ and ${\bf K}(\mu_{p^m})$ is
linearly disjoint over ${\bf K}$ with the compositum of the others.
Let $I$ be a subset of $\{1,...,n\}$ and $J$ a subset of $\{1,...,r\}$.
Then there exists a Cebotarev set $X$ of primes $w \not \in S$ such that
\newline
\noindent
1) $w$ is $\rho_m$-nice.
\newline
\noindent
2) $f_i |_{G_w} \neq 0$  for $i \in I$ and
 $f_i |_{G_w} = 0$ for $i \in \{1,\dots,n\} \backslash I$.
\newline
\noindent
3) $\phi_j|_{G_w} \neq 0$ for $j \in J$ and
$\phi_j|_{G_w} =0$ for $j \in \{1,\dots,r\} \backslash J$.
\end{fact}
Proof: This is a minor variant lemma $8$ of [Kh-Ram].
\hfill $\square$
\vskip1em

Let $\rhob$ be as given and $S$ be a set
containing $p$ and the ramified primes 
of $\rhob$.
We first enlarge $S$ so that global obstructions to deformation
questions can be locally detected. Recall that by definition
the symbol $\sha^i_{S }(\ad)$ is the kernel of the map
$H^i(G_{S},\ad)\rightarrow \oplus_{v \in S } H^i(G_v,\ad)$.
\begin{lemma}\label{Sha2}
Let $\rhob$ and $S$ be as above. There exists a finite
set $T$ of nice primes such that 
$\sha^1_{S \cup T}(\ad)$ and
$\sha^2_{S \cup T}(\ad)$ are trivial.
\end{lemma}
Proof: Let $M=\ad$ (resp. $M=\adst$)  and
let  ${\bf K}={\mathbb Q}(\ad,\mu_p)$.
Let $G=Gal({\bf K}/{\mathbb Q})$. That the image of $\rhob$
contains $SL_2({\mathbb Z}/p{\mathbb Z})$ implies
$H^1(G,M)=0$ in both cases. In each 
case we show that $\sha^1_{S \cup Y}(M) \subset \sha^1_S(M)$
for {\em any} set of primes $Y$. Indeed, for $\alpha
\in \sha^1_{S \cup Y}(M)$, we see using that
$H^1(G,M)=0$ and  Fact~\ref{disjointness}
that $\alpha$ cuts out a nontrivial 
extension ${\bf K}_{\alpha}$ of ${\bf K}$
such that $Gal({\bf K}_{\alpha}/{\bf K}) \simeq \ad$ (resp. $\adst$)
as a ${\mathbb Z}/p{\mathbb Z}[Gal({\bf K}/{\mathbb Q})]$-module.
Since $\alpha \in \sha^1_{S \cup Y}(M)$ we see 
the extension ${\bf K}_{\alpha}/{\bf K}$ is trivial at all places
of $S \cup Y$, and is therefore unramified  at all primes of $Y$.
Thus $\alpha$ inflates from an element of $H^1(G_S,M)$ that is trivial
at all places of $S$, so $\alpha \in \sha^1_S(\ad)$.

By global duality
$\sha^2_{S}(\ad)$ and $\sha^1_S(\adst)$ are dual.
Choose bases $\{f_1,...,f_m\}$ and $\{ \phi_1,...,\phi_r\}$
of $\sha^1_S(\ad)$ and $\sha^1_S(\adst)$ and 
sets of nice primes $\{ a_1,...,a_m\}$ and
$\{b_1,...,b_r\}$ are as in Fact~\ref{disjointness}
such that $f_i(\sigma_{a_i}) \neq 0$, $f_s(\sigma_{a_i})=0$ for $s\neq i$,
and $\phi_j(\sigma_{b_j}) \neq 0$, 
$\phi_t(\sigma_{b_j})=0$ for $t\neq j$.
\hfill $\square$
\vskip1em\noindent
Henceforth we enlarge $S$ as in the lemma.

\begin{lemma}\label{increase}
Let $\rhob$ be as usual, and suppose $S$ is such that $\sha^1_S(\ad)$
and $\sha^2_S(\ad)$ are trivial. 
For  $r \not \in S$ 
the inflation maps
\begin{equation}H^1(G_S,\ad) \rightarrow H^1(G_{S \cup \{r\}},\ad)
\end{equation}
and
\begin{equation}
H^1(G_S,\adst) \rightarrow H^1(G_{S \cup \{r\}},\adst)\end{equation}
have
cokernels of dimension $\dim(H^2(G_r,\ad))$ 
and $\dim(H^2(G_r,\adst))$
respectively and are injective. 
If $r$ is nice both cokernels are one dimensional.
\end{lemma}
Proof: The inflation maps are 
necessarily injective. Since $\sha^2_S(\ad)$
is trivial, by global duality we have $\sha^1_S(\adst)=0$. 
For all $v\in S \cup \{r\}$ put ${\cal L}_v = H^1(G_v,\ad)$ 
in Fact~\ref{wiles1.6}. 
Then $H^1_{\cal L}(G_S,\ad) = H^1(G_S,\ad)$ 
and
$H^1_{ {\cal L}^{\perp} }(G_S,\adst) =\sha^1_S(\adst)=0$.
Adding in the prime $r$ changes the total contribution to
the right hand side of Equation (\ref{fact4}) for $\ad$ by
$\dim(H^1(G_r,\ad))-\dim(H^0(G_r,\ad))$ which by the local 
Euler-Poincar\'{e} characteristic is just $\dim(H^2(G_r,\ad)$.
Since {\em adding} more primes to the ramified set 
cannot cause  the groups $\sha^1$ to increase in size we are done
for $\ad$. The proof for $\adst$ follows similarly
using the local conditions ${\cal M}_v = H^1(G_v,\adst)$.
If $r$ is nice Lemma~\ref{nice} completes the proof.
\hfill $\square$.

\vskip1em\noindent
\begin{lemma}\label{unram} 
Let $\rho_m$ be a deformation of $\rhob$ 
to $GL_2({\mathbb Z}/p^m{\mathbb Z})$
unramified outside a set $S$ and assume
$\sha^1_S(\ad)$ and $\sha^2_S(\ad)$ are trivial.
Let $R$ be any finite collection
of unramified primes of $\rhob$ disjoint from $S$.
Then there is a finite set  $Q=\{q_1,...,q_n\}$
of $\rho_m$-nice primes disjoint from $R \cup S$
such that
the maps

\begin{equation}\label{thirteen}
H^1(G_{S \cup R\cup Q},\adst) \rightarrow \oplus_{v \in Q}
H^1(G_v,\adst),\end{equation}

\begin{equation}\label{twelve}
H^1(G_{S \cup R\cup Q},\ad) \rightarrow \oplus_{v \in S \cup R}
H^1(G_v,\ad)\end{equation}
and
\begin{equation}H^1(G_{S \cup Q},\ad) \rightarrow
\left(\oplus_{v \in S}  H^1(G_v,\ad) \right)
\oplus \left(  \oplus_{r\in R} H^1_{nr}(G_r,\ad)
\right)\end{equation} are isomorphisms.
\end{lemma}
Proof: 
Let $\{\phi_1,...,\phi_n\}$  be a basis of $H^1(G_{S \cup R},\adst)$
and, using Fact~\ref{disjointness}
for $i=1,...,n$, choose $q_i$ that are $\rho_m$-nice 
and satisfy the conditions $\phi_i(\sigma_{q_i}) \neq 0$
and $j \neq i$ implies $\phi_j(\sigma_{q_i})=0$.
Put $Q=\{q_1,...,q_n\}$
and ${\cal L}_{q_i} = H^1(G_{q_i},\ad)$.
For $v \in S \cup R$ put ${\cal L}_v=0$.
Consider the dual Selmer map
$$H^1(G_{S \cup R\cup Q},\adst) \rightarrow$$
\begin{equation}\left(\oplus_{v \in S \cup R}  \frac{H^1(G_v,\adst)}{H^1(G_v,\adst)} \right)
\oplus
\left(\oplus^n_{i=1 }  \frac{H^1(G_{q_i},\adst)}{0}\right)\end{equation}
for the set $S\cup R \cup Q$ and the conditions ${\cal L}^{\perp}_v$.
Any  $\phi \in H^1_{ {\cal L}^{\perp} }(G_{S \cup R\cup Q},\adst)$
is clearly unramified at all primes of $Q$
and so inflates from $H^1(G_{S \cup R},\adst)$. 
Then $\phi$ is a linear combination
of the $\phi_i$ and necessarily non-trivial at some $q_i$. 
This contradiction
shows $\phi=0$  so 
$H^1_{ {\cal L}^{\perp} }(G_{S \cup Q \cup R},\adst)=0$.

We now show $H^1_{ {\cal L} }(G_{S \cup Q \cup R},\ad)=0$.
By assumption $H^1_{ {\cal L}}(G_S,\ad)=0$.
We use Fact~\ref{wiles1.6} with $M=\ad$ and the sets
$S \cup R$ and $S\cup R \cup Q$ respectively.
Then the change in the right hand side of Equation (\ref{fact4})
is 
\begin{equation}\sum_{q_i \in Q} 
\left(\dim(H^1(G_{q_i},\ad)) -\dim(H^0(G_{q_i},\ad))\right)
=\sum_{q_i \in Q} (2-1) =n.\end{equation}
The left hand side of Equation (\ref{fact4})
for the set $S \cup R$ is just $0-n$ so
for the set $S \cup R \cup Q$ the left hand side is $0-n+n=0$.
Since we already have
$H^1_{ {\cal L}^{\perp} }(G_{S \cup Q \cup R},\adst)=0$
it follows that $H^1_{ {\cal L} }(G_{S \cup Q \cup R},\ad)=0$
as well.

Using Fact~\ref{poitoutate} twice (once with local conditions
${\cal L}_v$ for $\ad$ and once with local conditions ${\cal L}^{\perp}_v$
for $\adst$) we get that
\begin{equation}H^1(G_{S \cup R\cup Q},\ad) \rightarrow \oplus_{v \in S \cup R}
H^1(G_v,\ad)\end{equation}
and 
\begin{equation}H^1(G_{S \cup R\cup Q},\adst) \rightarrow \oplus_{v \in Q}
H^1(G_v,\adst)\end{equation}
are isomorphisms.

Consider the inverse image in 
Equation (\ref{twelve}) 
of 
\begin{equation}\oplus_{v \in S}  H^1(G_v,\ad)
\oplus  \left(\oplus_{r\in R} H^1_{nr}(G_r,\ad)
\right).\end{equation} This inverse image clearly
lies in $H^1(G_{S \cup Q},\ad)$ which by
Lemma~\ref{increase} is of codimension 
$\oplus_{r \in R}\dim(H^2(G_r,\ad)$ in
$H^1(G_{S \cup Q \cup R},\ad)$. By the local 
Euler-Poincar\'{e} characteristic
we know $\oplus_{r\in R}H^1_{nr}(G_r,\ad)$ is of codimension 
$\oplus_{r \in R}\dim(H^2(G_r,\ad)$ in 
$\oplus_{r \in R}H^1(G_r,\ad)$ so the inverse image is exactly
$H^1(G_{S \cup Q},\ad)$. 
\hfill $\square$.
\vskip1em\noindent
Remark: The utility of Lemma~\ref{unram} is that it provides
us with global cohomology classes that will serve to `correct'
the local deformation problems at primes of $S$ and remove obstructions
to lifting at these primes, 
and to also arrange the characteristic 
polynomials at the unramified
primes of $R$ to be what we want. 
All this is done at the cost of possibly
introducing local obstructions to lifting at the nice primes of $Q$. 
These last obstructions will be removed by introducing for each
$q_i \in Q$ up to two nice primes so that a cohomology class
ramified at these new primes will correct things at $q_i$
{\em without} changing anything at primes in 
$S \cup R \cup Q \backslash \{q_i\}$ or introducing
obstructions at the new primes.

\begin{prop}\label{cheb}
Let $\rho_m$ be a deformation of $\rhob$ to $GL_2({\mathbb Z}/p^m{\mathbb Z})$
unramified outside a set $S$ and assume
$\sha^1_S(\ad)$ and $\sha^2_S(\ad)$ are trivial.
Let $S,R$, and $Q$ as in Lemma~\ref{unram}. 
We write $Q=\{q_1,...,q_n\}$. 
Let $A$ be any finite set of primes disjoint
from $S \cup R \cup Q$.
Fix $k$ between $1$ and $n$.
There exists
a Cebotarev set $T_{k}$ of primes $t_{k}$ such that
\begin{itemize}
  \item all $t_{k}\in T_k$ are $\rho_m$-nice
  \item for any $t_k \in T_k$, the kernel of the map 
        \begin{equation}H^1(G_{S \cup Q \cup \{t_{k}\}},\ad) \rightarrow
        \oplus_{v \in S } H^1(G_v,\ad) \oplus_{v \in R} H^1_{nr}(G_v,\ad)\end{equation}
        is one dimensional, spanned by $f_{t_k}$.
  \item $f_{t_k}|_{G_v}=0$ for all $v \in S\cup R \cup Q\cup A \backslash
        \{q_k\}$ and $f_{{t_k}}$ is unramified
        at $G_{q_k}$ with $f_{t_k}|_{G_{q_{k}}} \neq 0$ 
\end{itemize}
\end{prop}
Proof: 
Equation (\ref{thirteen}) identifies $H^1(G_{S \cup Q \cup R},\adst)$
with $\oplus^n_{i=1} H^1(G_{q_i},\adst)$. Since each summand
on the right hand side of Equation (\ref{thirteen})
is by Lemma~\ref{nice}
two dimensional, say with basis 
$\{\phi_{i1},\phi_{i2} \}$ where $\phi_{i1}$ spans
$H^1_{nr}(G_{q_i},\adst)$, we can abuse notation and consider
the union of these sets as $i$ runs from $1$ to $n$ to be a basis
of the left hand side. Thus
$H^1(G_{S \cup R \cup Q \cup A},\adst)$ has dimension $2n+d$ for some $d$.

We will use Fact~\ref{wiles1.6} with the set $S \cup R \cup Q \cup A$
and 
the local conditions
\begin{equation}{\cal L}_{q_k} =
H^1_{nr}(G_{q_{k}},\ad),\,\,\,\,{\cal L}_v =0 \,\,\,\,otherwise.\end{equation}
It is well known under the local pairing that ${\cal L}^{\perp}_{q_k} 
=H^1_{nr}(G_{k_k},\adst)$. For $v \neq q_k$ the spaces
${\cal L}^{\perp}_v$ are obvious.

The isomorphism of Equation (\ref{twelve})
implies
$H^1_{\cal L}(G_{S \cup R \cup Q \cup A},\ad)=0$. By definition 
$H^1_{{\cal L}^{\perp}}(G_{S \cup R \cup Q \cup A},\adst)$ 
is the kernel of the map
\begin{equation}H^1(G_{S \cup R \cup Q \cup A},\adst) \rightarrow 
\frac{ H^1(G_{q_k},\adst) }{H^1_{nr}(G_{q_k},\adst)}\end{equation}
whose target is one dimensional by Lemma~\ref{nice}.
Since $\phi_{k2}$ does not go to $0$ the map above is surjective
so $H^1_{{\cal L}^{\perp}}(G_{S \cup R \cup Q \cup A},\adst)$
has dimension $2n+d-1$. It has the basis
$\{ \phi_{11},\phi_{12},...,\phi_{n1},\phi_{n2},\psi_1,...\psi_d\}
\backslash \{ \phi_{k2}\}$.

Using Fact~\ref{disjointness}
let $T_k$ be the Cebotarev set of $\rho_m$-nice primes $t_k$
such that all the basis elements of 
$H^1_{{\cal L}^{\perp}}(G_{S \cup R \cup Q \cup A},\adst)$
are trivial at $\sigma_{t_k}$ {\em and}
$\phi_{k2}(\sigma_{t_k}) \neq 0$.
We extend out local conditions ${\cal L}_v$ to $t_k$
by 
putting ${\cal L}_{t_k}=H^1(G_{t_k},\ad)$ so
${\cal L}^{\perp}_{t_k}=0$.

We know by Fact~\ref{wiles1.6} that
\begin{equation}\dim(H^1_{\cal L}(G_{S\cup R \cup Q \cup A},\ad)) 
-\dim(H^1_{{\cal L}^{\perp}}(G_{S\cup R \cup Q \cup A},\adst)) =-2n-d+1\end{equation}
so adding in the prime $t_{k}$
gives a contribution to the right hand side
of Equation (\ref{fact4}) of $2-1=1$ so
$$\dim(H^1_{\cal L}(G_{S\cup R \cup Q \cup A\cup \{t_k\}},\ad)) 
-\dim(H^1_{{\cal L}^{\perp}}(G_{S\cup R \cup Q \cup A\cup \{t_k\}},\adst))$$
\begin{equation}=-2n-d+2.\end{equation}

Since ${\cal L}^{\perp}_{t_{k}}=0$, we see any element in 
$H^1_{{\cal L}^{\perp}}(G_{S\cup R \cup Q \cup A\cup \{t_k\}},\adst)$
is unramified at $t_{k}$ and thus inflates
from $H^1_{{\cal L}^{\perp}}(G_{S \cup R \cup Q \cup A},\adst)$.
This last group has 
basis $\{ \phi_{11},\phi_{12},...,\phi_{n1},\phi_{n2},\psi_1,...\psi_d\}
\backslash \{ \phi_{k2}\}$ and
these are all  trivial at $t_{k}\in T_k$.
Thus $H^1_{{\cal L}^{\perp}}(G_{S \cup R \cup Q \cup A\cup \{ t_{k} \} },\adst)$
has dimension $2n+d-1$ as well and so
$H^1_{{\cal L}}(G_{S \cup R \cup Q \cup A \cup \{ t_{k} \} },\ad)$
has dimension $1$. Let $f_{t_{k}}$ span this Selmer group.
Note $f_{t_k}$ is by definition unramified at $q_k$ and
clearly $f_{t_k}$ is trivial at all primes in 
$R \cup A $
and therefore it inflates from $H^1(G_{S \cup Q \cup \{t_k\}},\ad)$.

It remains to show 
$f_{t_k}|_{G_{q_k}} \neq 0$. We show 
$f_{t_k}(\sigma_{q_k}) \neq 0$. To do this we introduce
a slightly different set of local conditions. 
Put ${\cal M}_v=0$ for all $v \in S \cup R \cup Q \cup A$.
Then $H^1_{ {\cal M}^{\perp} }(G_{S \cup R \cup Q \cup A},\adst)$
is just  $H^1(G_{S \cup R \cup Q \cup A},\adst)$ and has basis
$\{\phi_{11},\phi_{12},...\phi_{n1},\phi_{n2},
\psi_1,...,\psi_d \}$ with $2n+d$ elements.
We already observed 
$H^1_{\cal L}(G_{S \cup R \cup Q \cup A},\ad)=0$ so we see
$H^1_{\cal M}(G_{S \cup R \cup Q \cup A},\ad)=0$.
Thus
\begin{equation}\dim(H^1_{\cal M}(G_{S\cup R \cup Q \cup A},\ad))
-\dim(H^1_{{\cal M}^{\perp}}(G_{S\cup R \cup Q \cup A},\adst)) =-2n-d.\end{equation}
Adding in the prime $t_{k}$ with 
${\cal M}_{t_k}={\cal L}_{t_k} =H^1(G_{t_k},\ad)$
gives a new local contribution of $2-1=1$ 
in the right hand side of Equation (\ref{fact4})
so
$$\dim(H^1_{\cal M}(G_{S\cup R \cup Q \cup A\cup \{t_k\}},\ad))
-\dim(H^1_{{\cal M}^{\perp}}(G_{S\cup R \cup Q \cup A\cup \{t_k\}},\adst))$$
\begin{equation}=-2n-d+1.\end{equation}
As ${\cal M}^{\perp}_{t_k}=0$, we see any element
of $H^1_{{\cal M}^{\perp}}(G_{S\cup R \cup Q \cup A\cup \{t_k\}},\adst)$
is unramified at $t_k$ and inflates from
$H^1_{{\cal M}^{\perp}}(G_{S\cup R \cup Q \cup A},\adst)$
which is spanned by
$\{\phi_{11},\phi_{12},...,\phi_{n1},\phi_{n2},
\psi_1,...,\psi_d \}$.
By our choice of $t_k \in T_k$ we see all of these basis elements
{\em except} $\phi_{k2}$ are in 
$H^1_{{\cal M}^{\perp}}(G_{S\cup R \cup Q \cup A\cup \{t_k\}},\adst)$
so this dual Selmer group has dimension $2n+d-1.$
This implies
$H^1_{\cal M}(G_{S\cup R \cup Q \cup A\cup \{t_k\}},\ad)$ is trivial
and therefore does {\em not} contain $f_{t_k}$. That is
$f_{t_k} |_{G_{q_k}} \neq 0$ and we are done.
\hfill $\square$.
\vskip1em
We will choose up to two primes from each $T_k$ 
and the  corresponding cohomology classes will 
be used
to remove any obstructions to lifting at $q_k$
{\em without}
introducing obstructions at the prime(s) chosen
from $T_k$. We currently do not know how to
do this with just one prime from each $T_k$. 
It would be of interest to find a method
that works for one prime.
Proposition~\ref{cheb} also guarantees that
the cohomology class we use will {\em not} change anything
at primes in $S\cup R$ or at any of the primes $q_i$ for $i \neq k$.
Finally, we need to guarantee, that as $i$ runs
from $1$ to $n$, using a prime(s) $t_i \in T_i$ and its cohomology
class will not change things locally at 
the union over $j<i$ of the prime(s) in $T_j$.
Allowing Proposition~\ref{polarisation} below, one can proceed
to the next section.

\begin{prop}\label{polarisation}
Let $T_k$ be as in Proposition~\ref{cheb}.
Then there is a set $\tilde{T}_k \subset T_k$ of one or two primes such that
\begin{itemize}

  \item There is a linear combination $f_k$ of the elements $f_{t_k}$ for
        $t_k \in \tilde{T}_k$ such that $f_k(\sigma_{t_k})=0$ for
        all $t_k \in \tilde{T}_k$ and $f_k|_{G_{q_k}} \neq 0$
  \item $j<k$ implies that for 
        $t_j \in \tilde{T}_j$ we have $f_{t_k}(\sigma_{t_j})=0$
  \item $f_k|_{G_v}=0$ for all $v \in S \cup R \cup Q \backslash \{q_k\}$
\end{itemize}
\end{prop}
Proof: 
We refer the reader to Remark $1$) 
following Lemma~\ref{nice} for the definition
of $f_{t_k}(\sigma_{t_k})$.
We induct.
Suppose $\tilde{T}_i$ and $f_i$ have been chosen for $i=1,...,k-1$.
We show how to choose $\tilde{T}_k$ and $f_k$. 
Put $A_k = \cup^{k-1}_{i=1} \tilde{T}_i$. Applying Proposition~\ref{cheb}
with $A_k$  here playing the role of $A$ there
we get a Cebotarev
set $T_k$ satisfying the last two
required properties. It remains to ensure the first.

Consider the map
$$H^1(G_{S \cup R \cup Q \cup A_k},\adst) \rightarrow $$
\begin{equation}\label{withouttk} 
\left(\oplus_{v \in S \cup R \cup Q \cup A_k \backslash \{q_k\} } 
\frac{H^1(G_v,\adst)}{H^1(G_v,\adst)}\right)
\oplus \frac{H^1(G_{q_k},\adst}{H^1_{nr}(G_{q_k},\adst)}.\end{equation}
Its kernel, $H^1_{ {\cal L}^{\perp}}(G_{S \cup R \cup Q \cup A_k},\adst)$,
has dimension $2n-1+d_k$ for some $d_k \geq 0$
and is spanned by all the $\phi_{ij}$ and $\psi_i$
{\em except} $\phi_{k2}$. In particular the map is surjective.
Now consider 
$$H^1(G_{S \cup R \cup Q \cup A_k \cup \{t_k\}},\adst) 
\rightarrow $$ 
\begin{equation}\label{withtk}
\left(\oplus_{v \in S \cup R \cup Q \cup A_k \backslash \{q_k\} }  
\frac{H^1(G_v,\adst)}{H^1(G_v,\adst)}\right)
\oplus \frac{H^1(G_{q_k},\adst}{H^1_{nr}(G_{q_k},\adst)}.\end{equation}
Lemma~\ref{increase} implies
$H^1(G_{S \cup R \cup Q \cup A_k \cup \{t_k\}},\adst)$
has dimension one bigger than
$H^1(G_{S \cup R \cup Q \cup A_k },\adst)$,
so the kernel of Equation (\ref{withtk})
is one dimension bigger than the kernel of 
Equation (\ref{withouttk}). Let $\phi_{t_k}$
be {\em any} element in the kernel of this second map that is not
in the kernel of the first map. Note $\phi_{t_k}$ is necessarily
ramified at $t_k$.
We will need $\phi_{t_k}$ later. It is not
well defined, but as it turns out
this ambiguity will not matter to us.

Suppose there is a prime $t_k \in T_k$ such that
$f_{t_k}(\sigma_{t_k})=0$. Since $t_k$ is $\rho_m$-nice,
adjusting $\rho_m$ by a multiple of $f_{t_k}$ will keep
$\rho_m|_{G_{t_k}}$ unobstructed as 
in Definition~\ref{unobstructed}. In this case
we can take $\tilde{T}_k = \{ t_k\}$.
Henceforth we assume that for all $t_k \in T_k$ that
$f_{t_k}(\sigma_{t_k}) \neq 0$.
Also, since by Proposition~\ref{cheb} we know
$f_{t_k}(\sigma_{q_k}) \neq 0$, by replacing
$f_{t_k}$ by an appropriate multiple we may assume
$f_{t_k}(\sigma_{q_k})=1$ for all $t_k \in T_k$.

It turns out we need to find $t_{k1}$ and $t_{k2}$ in $T_k$ such that
the $2 \times 2$ matrix $(f_{t_i}(\sigma_{t_j}))_{1 \leq i,j \leq 2}$
has determinant
$0$ {\em and} unequal rows.
The matrix, by assumption has non-zero diagonal entries.
Suppose the matrix is as below
\vskip3em
\hspace{1.5in}
\begin{tabular}{c|c|c|c}
           &     $\sigma_{t_{k1}}$   & $\sigma_{t_{k2}}$   &$\sigma_{q_k}$ \\
\hline
$f_{t_{k1}}$  &     $a$            &     $b$        &   $1$\\
\hline
$f_{t_{k2}}$  &     $c$            &     $d$        &   $1$  \\
\end{tabular}
\vskip3em

Since $f_{t_{ki}}|_{G_v}=0$ for all $v \in S \cup
Q\cup R \cup A_k \backslash \{q_k\},$
adjusting our mod $p^m$
representation by a linear combination of
$f_{t_{k1}}$ and $f_{t_{k2}}$ could cause new local
obstructions only at  $t_{k1}$ and $t_{k2}$.
We want
$f_k:=\alpha_1f_{t_{k1}}+
\alpha_2f_{t_{k2}}$ for  $\alpha_i \in {\mathbb Z}/p{\mathbb Z}$
to unobstruct the $q_k$-local deformation problem.
We require  $f_k(\sigma_{t_{ki}})=0$
for $i=1,2$, and for any $\beta \in {\mathbb Z}/p{\mathbb Z}$, we must
be able to solve $(\alpha_1f_{t_{k1}}+\alpha_2f_{t_{k2}})(\sigma_{q_k})=\beta$.
Showing that $\alpha_1$ and $\alpha_2$ exist as required
is equivalent to guaranteeing the conditions on the matrix described above.

Recall we are assuming that for $t_k \in T_k$, 
$f_{t_k}(\sigma_{t_k}) \neq 0$. Let $y$ be the
value of $f_{t_k}(\sigma_{t_k}) $ that occurs most often,
that is with maximal upper density.
Let
\begin{equation}
Y_k=\{t_k \in T_k | f_{t_k}(\sigma_{t_k})=y\}.
\end{equation} 
Then $Y_k$ may not have a density,
but it has a positive upper density.

By the proof of Lemma~\ref{nice}
\begin{equation}H^1(G_{t_k},\ad) \simeq H^1(G_{t_k},{\mathbb Z}/p{\mathbb Z})
\oplus H^1(G_{t_k},\mu_p)\end{equation}
so  we can consider $f_{t_k,{\mathbb Z}/p{\mathbb Z}}$ and
$f_{t_k,\mu_p}$,
the projections of $f_{t_k}|_{G_{t_k}}$ to the direct summands. 
Since  $f_{t_k}$ is necessarily
ramified at $t_k$
we see $f_{t_k,\mu_p} \neq 0$. Using Lemma~\ref{nice}
again we can decompose 
elements of $H^1(G_{t_k},\adst)$ similarly.

Recall that for any prime $r$ that
$H^2(G_r,\mu_l)$ is the $l$-torsion in ${\mathbb Q}/{\mathbb Z}$,
the Brauer group of ${\mathbb Q}_r$.
For  any
nice prime $q$, define $g_q \in H^1(G_q,{\mathbb Z}/p{\mathbb Z})$
by $g_q(\sigma_q)=1$.
For all $t_k \in Y_k$ consider the necessarily nonzero values 
$f_{t_k,\mu_p} \cup g_{t_k}$ in
${\mathbb Q}/{\mathbb Z}$. Let $z$
be the value that occurs most often.
Put \begin{equation}Z_k=  
\{ t_k \in Y_k | f_{t_k,\mu_p} \cup g_{t_k}=z \}.\end{equation}
Then $Z_k$ has positive upper density.
Note that $Z_k$ depends only on $T_k$ and the prime $q_k$
and that $ Z_k \subset Y_k \subset T_k$.

Choose any $t_{k1} \in Z_k$.
Recall $\phi_{t_{k1}}$ is defined (with $t_{k1}$ playing
the role of $t_k$) as an element in the kernel of Equation (\ref{withtk})
that is {\em not} in the kernel of Equation (\ref{withouttk}).
We will try to choose $t_{k2} \in Z_k$
so our $2 \times 2$ matrix
$(f_{t_{ki}}(\sigma_{t_{kj}}))_{ 1 \leq i,j \leq 2}$
has the desired properties. As $t_{k1},t_{k2} \in Z_k$, both diagonal
entries
will both be $y$.
Choosing $f_{t_{k1}}(\sigma_{t_{k2}})$ to be what we want (say $x \neq 0,y$)
is simply a Cebotarev
condition on $t_{k2}$
independent of those that determine $T_k$.
Choosing $f_{t_{k2}}(\sigma_{t_{k1}})$ as we want
($y^2/x$ in this case) involves 
invoking the global reciprocity law to make the choice
a Cebotarev condition.

Consider the element of $H^2(\gal, \mu_p)$ given
by $\phi_{t_{k1}} \cup f_{t_{k2}}$. This class is unramified
outside $S \cup Q \cup R \cup A_k \cup \{t_{k1},t_{k2}\}$.
Recall there was ambiguity in the definition of $\phi_{t_{k1}}$
(see Equation (\ref{withtk})),
namely we have no control over its local behavior at 
$v \in S \cup R \cup Q \cup A_k \backslash \{q_k\}$.
But these are precisely the places at which $f_{t_{k2}}$ is
trivial and since we will be summing local invariants,
the ambiguity is irrelevant.
Also, $\phi_{t_{k1}}$ and $f_{t_{k2}}$ are both unramified
at $q_k$ and thus at $q_k$ their cup product is zero as well.
As the sum of the local invariants is zero we have
\begin{equation}\label{invariant} 
inv_{t_{k1}}(\phi_{t_{k1}} \cup f_{t_{k2}}) = 
- inv_{t_{k2}}(\phi_{t_{k1}} \cup f_{t_{k2}}).\end{equation}
Consider the left hand side and recall we are supposing
we have chosen $t_{k1}$ and are trying to choose $t_{k2}$.
Then $\phi_{t_{k1}}$ is fixed and the left hand side
depends entirely on $f_{t_{k2}}(\sigma_{t_{k1}})$. Thus 
choosing the left hand side to be whatever we want is 
equivalent to
choosing $f_{t_{k2}}(\sigma_{t_{k1}})$ to be whatever we want.

Now consider the right hand side of Equation (\ref{invariant}). 
Recall $g_{t_{k1}} \in H^1(G_{t_{k1}},{\mathbb Z}/p{\mathbb Z})$ 
is normalised so $g_{t_{k1}}(\sigma_{t_{k1}})=1$,
$t_{k2} \in Z_k$,  $\phi_{t_{k1}}$ is unramified at $t_{k2}$
and $f_{t_{k2}}$ is ramified at $t_{k2}$. Thus we have
$$ inv_{t_{k2}}(\phi_{t_{k1}} \cup f_{t_{k2}})
=inv_{t_{k2}}(\phi_{t_{k1},{\mathbb Z}/p} \cup 
f_{t_{k2},\mu_p})=$$
\begin{equation}
 \phi_{t_{k1}}(\sigma_{t_{k2}}) \cdot inv_{t_{k2}}
(g_{t_{k1}} \cup f_{t_{k2},\mu_p})
=\phi_{t_{k1}}(\sigma_{t_{k2}})z.
\end{equation}
Thus the left hand side of Equation (\ref{invariant}) depends
only on $\phi_{t_{k1}}(\sigma_{t_{k2}})$.
So 
choosing $f_{t_{k2}}(\sigma_{t_{k1}})$ to be
whatever we like 
is equivalent to choosing  $\phi_{t_{k1}}(\sigma_{t_{k2}})$
to be whatever we like.

So, if we can choose $t_{k2} \in Z_k$ such that
$f_{t_{k1}}(\sigma_{t_{k2}})$ and $\phi_{t_{k1}}(\sigma_{t_{k2}})$ are whatever
we wish, we'll be able to choose $f_{t_{k1}}(\sigma_{t_{k2}})$ and
$f_{t_{k2}}(\sigma_{t_{k1}})$ to be a non-zero $x\neq y$ and $y^2/x$ respectively
and we'll be done. (By Fact~\ref{disjointness}, $f_{t_{k1}}$ and
$\phi_{t_{k1}}$ give  independent Cebotarev conditions.)

Suppose, having chosen $t_{k1}$, we can't do this. Then
the set  $Z_k \backslash \{ t_{k1}  \}$ 
lies in Cebotarev classes that are complementary
to the Cebotarev conditions on $\sigma_{t_{k2}}$ imposed
by choosing $f_{t_{k1}}(\sigma_{t_{k2}})=x$ where $x \neq 0,y$
and choosing $\phi_{t_{k1}}(\sigma_{t_{k2}})$ to be whatever 
forces $f_{t_{k2}}(\sigma_{t_{k1}})=y^2/x$.
Note that if the set $T_k$ has
density $D$, then these complementary Cebotarev classes
form a set of density $D\gamma$ where $\gamma = 1-\frac{p-2}{p^2}$.
(The actual value of $\gamma$ is not important; all that matters
is that $\gamma <1$.)

Now replace $t_{k1}$ by a sequence of different primes $l\in Z_k$, and
assume they also allow no valid choice for the second prime.
Then we see that $Z_k  \backslash \{ l \}$ 
also lies in the complimentary Cebotarev classes
associated to $f_l$ and $\phi_l$. But these classes, for varying $l$, are
all independent of one another ($\phi_l$ and $f_l$ being ramified at $l$), 
so imposing $n$ such conditions, the density of the complimentary classes
is $D\gamma^n$. Thus we have that $Z_k \backslash \{ l_1,...,l_n \}$
is contained in a set of density $D\gamma^n$.
This holds for all positive $n$.
Letting $n$ get arbitrarily large
we get that
$Z_k $
is contained in a set of arbitrarily 
small density, so $Z_k$ has upper density
$0$, a contradiction.

We can choose primes $\{t_{k1}, t_{k2} \}$ so that our matrix has the desired
properties. 
Thus there is an
$f_k:=\alpha_1f_{t_{k1}}+\alpha_2f_{t_{k2}}$ that does
what we want.
The induction is complete, and the proposition is proved.
\hfill $\square$

\subsection{Application I}

In this section we prove the main theorem.

Let $\rhob$ be given. 
Suppose $det\rhob = \chi^k$, where $\chi$ is the mod
$p$ cyclotomic character and $1 \leq k \leq p-1$. We fix the determinant
of our deformation to be $\epsilon^k$ where $\epsilon$
is the $p$-adic cyclotomic character.
Enlarge $S$ so that $\sha^1_S(\ad)$
and $\sha^2_S(\ad)$ are trivial as in Lemma~\ref{Sha2}.

Since $\sha^2_S(\ad)=0$ and  {\em all}
local deformation problems have no obstruction to
lifting to $GL_2({\mathbb Z}/p^2{\mathbb Z})$,
we choose $\rho_2:G_S \rightarrow GL_2({\mathbb Z}/p^2{\mathbb Z})$
a deformation of $\rhob$. Put $S_2=S$.
Once and for all, choose {\em any} deformation of $\rhob|_{G_v}$
to $GL_2({\mathbb Z}_p)$ for $v \in S_2$.
Let $R_2$ be all primes not in $S_2$ that are less than
$\displaystyle {(p^2/2)}^{\frac{2}{k}}$.
Once and for all, choose characteristic polynomials
in ${\mathbb Z}[x]$ for the primes of $R_2$. We {\em cannot}
guarantee these polynomials are pure of weight $k$.

Let $Q_2$ be the set from Proposition~\ref{unram}
for $S_2$ and $R_2$. There is an $f_2 \in 
H^1(G_{S_2 \cup Q_2},\ad)$ 
such that $(I+pf_2)\rho_2 |_{G_v}$
is the mod $p^2$ reduction of the preselected deformation
of $\rhob|_{G_v}$ to $GL_2({\mathbb Z}_p)$ for
all $v \in S_2$. Furthermore
the characteristic polynomials of
Frobenius of primes in $R_2$ in $(I+pf_2)\rho_2$
are as chosen. 
Using  Proposition~\ref{polarisation}, we can remove any
local obstructions to deforming to $GL_2({\mathbb Z}/p^3{\mathbb Z})$
at primes of $Q_2$, by allowing ramification at some more
$\rho_2$-nice primes.
The set of primes in question is $V_2:=\cup \tilde{T}_i$. 
There is a
$g_2 \in H^1(G_{S_2 \cup R_2 \cup Q_2 \cup V_2},\ad)$
that is trivial at primes of $S_2 \cup R_2$ such that
$(I+p(f_2+g_2))\rho_2$ 
\begin{itemize}
   \item is unramified outside $S_2 \cup Q_2 \cup V_2$
   \item is locally at $v \in S_2$ the mod $p^2$ reduction of
         the preselected deformations to ${\mathbb Z}_p$
   \item has the preselected characteristic polynomials of
         Frobenius for all primes in $\rho_2$
   \item is unobstructed at primes of $Q_2\cup V_2$
\end{itemize}
Since $\sha^2_{S_2 \cup Q_2 \cup V_2}=0$,
$(I+p(f_2+g_2))\rho_2$
can be deformed to $GL_2({\mathbb Z}/p^3{\mathbb Z})$.
Let $\rho_3$ be such a deformation of $(I+p(f_2+g_2))\rho_2$ unramified 
outside $S_2 \cup Q_2 \cup V_2$. Put $S_3=S_2 \cup Q_2 \cup V_2$
and let $R_3$ be the union of $R_2$ and all primes not in 
$S_3$ less than $\displaystyle {(p^3/2)}^{\frac{2}{k}}$

Once and for all, choose {\em any} deformation of 
$(I+p(f_2+g_2))\rho_2|_{G_v}$
to $GL_2({\mathbb Z}_p)$ for $v \in S_3 \backslash S_2$.
Once and for all, choose pure weight $k$ characteristic polynomials
in ${\mathbb Z}[x]$ for the primes of 
$R_3 \backslash R_2$. These polynomials
are necessarily of the form $x^2-a_rx+r^k$. 
So choosing the characteristic polynomials amounts to choosing
the $a_r$.
If the discriminant $a^2_r-4r^k<0$,
the roots will be pure of weight $k$. But $a_r$ is determined
mod $p^2$, so we may only alter $a_r$ by multiples
of $p^2$.
If $2r^{k/2}>p^2$, that is if 
$r> (p^2/2)^{2/k}$, we can choose $a_r$ such that
$a^2_r -4r^k<0$. All smaller primes are already in $R_2$,
so the characteristic polynomials of primes in $R_3 \backslash R_2$
can be chosen to be pure.

Now repeat the argument above, with $Q_3$ and $V_3$.
Let $\rho_4$ be a deformation of some 
$(I+p^2(f_3+g_3))\rho_3$, and continue. 
In general, the bound for primes of $R_n$ to guarantee purity
will be
$\displaystyle (p^n/2)^{\frac{2}{k}}$.

The induction is complete.

\vskip1em\noindent
{\bf Main Theorem } {\em Let $p \geq 5$ and
$\rhob:\gal \rightarrow GL_2({\mathbb Z}/p{\mathbb Z})$ be given such that
the image of $\rhob$ contains $SL_2({\mathbb Z}/p{\mathbb Z})$.
If $p=5$  assume further $\rhob$ is surjective.
Assume det$\rhob=\chi^k$ where $\chi$
is the mod $p$ cyclotomic character and $1 \leq k \leq p-1$.
Then there exists a deformation $\rho$ of $\rhob$ to ${\mathbb Z}_p$
such that
$\rho|_{G_p}$ is potentially semistable,
$\rho$ is unramified at a density $1$ set of primes $R$,
and for all but finitely many unramified primes $r$ 
the characteristic polynomial of Frobenius at $r$
is in ${\mathbb Z}[x]$  pure of weight $k$.}
%(in the sense of Definition \ref{purity} below).
%Furthermore, for the density $0$ infinite set of primes which we added
%to our ramification set, we may choose the local representation
%at any of these primes to be whatever we want, provided our choice
%is consistent with the local representation at the mod $p^m$ level
%these primes were introduced. In particular, if $\rhob|_{G_p}$
%is unramified, we may choose $\rho|_{G_p}$ to be unramified.
\vskip1em\noindent
Proof: All that remains to check is the density statement for ramified
primes. As $\rho$ is clearly irreducible, this follows
from [Kh-Raj]. 
\hfill $\square$
\vskip1em\noindent
Remark: Here, to appeal directly to [Kh-Raj], we are using the fact that we have forced ramification at all
auxiliary primes that are introduced in the lifting process. If we  do not use this we can
appeal to Proposition \ref{rub} instead. 
Alternatively, one can give a more self-contained proof
instead of using [Kh-Raj] simply because all auxiliary primes $q$ that are introduced
are such that the characteristic polynomial of Frobenius at these primes, in the sense of Definition \ref{charpoly},
is of a very constrained form.

\subsection{Application II}

\begin{cor} There exists a surjective map $\rho:\gal \rightarrow
SL_2({\mathbb Z}_7)$ unramified at $7$ and  the prime
$7$ is almost totally 
split in the field fixed by the kernel of $\rho$.
\end{cor}
Proof: In Section $8$ of [R1] it is recalled that
the polynomial
$f(x)=x^7-22x^6+141x^5-204x^4-428x^3+768x^2+320x-512$ of
Zeh-Marschke has splitting field with Galois group over ${\mathbb Q}$
equal to $PSL_2({\mathbb Z}/7{\mathbb Z})$. It is shown in [R1], following
Zeh-Marschke, that there is a Galois extension over ${\mathbb Q}$ 
with Galois group $SL_2({\mathbb Z}/7{\mathbb Z})$ containing this splitting field.
A factorisation of $f(x)$ mod $7$ (keeping in mind the discriminant
of $f(x)$ is $2^{50}19^4 367^2$) shows that Frobenius at $7$ has order
$3$ in the $PSL_2({\mathbb Z}/7{\mathbb Z})$  extension of  ${\mathbb Q}$
and thus has order
$3$ or $6$ in the $SL_2({\mathbb Z}/7{\mathbb Z})$ extension
of ${\mathbb Q}$. We do the deformations with determinant
fixed once and for all to be trivial and the local at $7$ deformation
to be the 
Teichm\"{u}ller lift. While $k=0$ is excluded in the Main Theorem, this is
only for choosing bounds for primes to obtain purity at a density one set
of primes. If one forgoes purity, the proof of the Main Theorem applies here.
\hfill $\square$
\vskip1em\noindent
Remark: It is a consequence of
the Fontaine-Mazur conjecture that for a number field ${\bf K}$
there  is no everywhere unramified $p$-adic analytic pro-$p$ extension
of ${\bf K}$. In the above example if ${\bf K}$ is the
$SL_2({\mathbb Z}/7{\mathbb Z})$ extension of ${\mathbb Q}$
and ${\bf L}$ is the field fixed by our surjective representation
to $SL_2({\mathbb Z}_7)$, then $L/K$ is an
example of a  $7$-adic analytic pro-$7$ extension in which
the primes of ${\bf K}$ above $7$ split completely. Of course ${\bf L}/{\bf K}$
is ramified at infinitely many primes other than $7$ and so does
not provide a counterexample to the Fontaine-Mazur  conjecture.
\vskip1em
Let $\rhob$ be odd and absolutely irreducible and modular
of square-free level and weight $2$.
Let $S$ be a set containing $p$ and the ramified
primes of $\rhob$. Let $T$ be a set of nice primes
such that the dual Selmer group
of [R3] for the set $S \cup T$
is trivial. Then the unique deformation of $\rhob$ to ${\mathbb Z}_p$
is denoted $\rho^{T-new}_{S \cup T}$ in [Kh-Ram]
and the corresponding ring $R^{T-new}_{S \cup T} \simeq {\mathbb Z}_p$ gives rise
to the unique `new at $T$' newform $g$ whose Galois
$p$-adic representation is congruent to $\rhob$ mod $p$.
\vskip1em\noindent
\begin{cor} With the set-up as above, there exists
a set $U$ consisting of at most two nice primes such that
$R^{T \cup U-new}_{S \cup T \cup U} \not \simeq {\mathbb Z}_p$.
There are at least two `new at $T \cup U$' newforms congruent to $f$.
\end{cor}
Proof: For any $m>1$ let $\rho_m$ be the mod $p^m$
redution of $\rho^{T-new}_{S \cup T}$.
We apply 
Proposition~\ref{polarisation}
with $S\cup T$ here playing the role of $S$ there and
$R$ being trivial. 
Observe Proposition~\ref{polarisation}
can be applied with $Q$ empty or with
$Q=\{q_1\}$ where $\{q_1\}$ is {\em any} $\rho_m$-nice prime.
Set $U=\{t_1,t_2\}$.

Then there is an $f=\alpha_1f_{t_1}+\alpha_2f_{t_2}$ 
where the $t_i$ are $\rho_m$-nice,
$f({\sigma_{t_i}})=0$ for $i=1,2$,
and $f|_{G_v}=0$ for $v \in S \cup T$.
Recall from [R3] the local conditions ${\cal N}_v$.
(See [R3] or the discussion prior to  Lemma~\ref{aux} for a detailed
discussion of the ${\cal N}_v$. Note that  $f({\sigma_{t_i}})=0$ for $i=1,2$
implies $f \in {\cal N}_{t_i}$ for $i=1,2$.)
Then the map
\begin{equation}\label{Nv}H^1(G_{S \cup T \cup U},\ad)\rightarrow
\oplus_{v \in S \cup T \cup U} 
\frac{ H^1(G_v,\ad) }{ {\cal N}_v}\end{equation}
has kernel containing $f$. (As in Proposition~\ref{polarisation}
we are assuming no one prime set $U$ works.)

Since the tangent space of
$R^{T \cup U-new}_{S \cup T \cup U}$ is the dimension
of the kernel of Equation (\ref{Nv}) and since the  $R=T$ 
theorem of Wiles and Taylor-Wiles applies
in this situation, we see
$R^{T \cup U-new}_{S \cup T \cup U} \not \simeq {\mathbb Z}_p$.
There are at least two
`new at $T \cup U$' new forms whose residual representations
are isomorphic to $\rhob$.
\hfill $\square$
\vskip1em\noindent
Remark: If the hypotheses of the last corollary are satisfied,
it would be of interest to describe those nice primes $t$
for which $f_t(\sigma_t)=0$. For instance, does the set
have positive density? Is it a Cebotarev set of some sort?

\subsection{Application III}
We refer the reader to Definitions~\ref{charpoly}-~\ref{compatible}.
We are going to prove that given a more or less arbitrary pair
of representations 
$\rhob_p:\gal \rightarrow GL_2({\mathbb Z}/p{\mathbb Z})$ and 
$\rhob_q:\gal \rightarrow GL_2({\mathbb Z}/q{\mathbb Z})$
for distinct primes $p$ and $q$ (with the residual 
images of $\rhob_p$ and $\rhob_q$ `large')
we can lift them to a density $1$ compatible pair of (infinitely ramified)
$p$-adic and $q$-adic representations as in Definition \ref{compatible}. This again answers a question 
of [Kh-Raj], and shows that unlike in the case of finitely ramified compatible pair of $p$-adic and $q$-adic representations where
one expects a motive to lurk behind them and connect them up,
an infinitely ramified density $1$ compatible pair of 
$p$-adic and $q$-adic representations can have little to do with each other 
besides the assumed relation of compatibility!
Note here that by compatible here we mean only a condition at a density one set of primes (that excludes
for instance all the primes ramified in either of the $p$-adic or $q$-adic representations). Thus we 
are not saying that we will get a strictly compatible lift where by strict one imposes conditions 
in particular at ramified primes outside $p$ and $q$ (in the first approximation), the condition being that
they come from the {\it same complex} representation of the
Weil-Deligne group at these primes, i.e., they have the same
Weil-Deligne parameter. It cannot be expected that arbitrary 
$\rhob_p$ and $\rhob_q$ as above have compatible `motivic lifts'
which in particular will be finitely ramified and strictly compatible. 
This is because
if such motivic lifts exists, this would imply by standard
expectations about motivic representations 
that strictly compatible lifts exist, and thus at whichever
primes $\rhob_p$ ramifies, there 
is a lift of $\rhob_q$ ramified at such a prime. This imposes 
conditions on $\rhob_q$ at such primes, 
conditions coming from structure of local Galois groups.
We do, however,  expect that finitely ramified pure rational 
$p$-adic representations
come from motives.

\begin{cor}\label{comp} Let $p > q \geq 5$ be primes and let
$\rhob_p:\gal \rightarrow GL_2({\mathbb Z}/p{\mathbb Z})$ and
$\rhob_q:\gal \rightarrow GL_2({\mathbb Z}/q{\mathbb Z})$ be 
representations whose mod $p$ and mod $q$ images contain 
$SL_2$. (If one of the primes
is $5$ we need the map to be surjective in that case.)
Say both are weight $k$ with $1\leq k \leq q-1$. 
Then there exist  potentially
semistable deformations of $\rhob_p$ and $\rhob_q$ respectively
$\rho_p:\gal \rightarrow GL_2({\mathbb Z}_p)$ and
$\rho_q:\gal \rightarrow GL_2({\mathbb Z}_q)$
such that for a set of primes $R$ of density one we have
that for $r \in R$ the characteristic polynomials
of $\sigma_r$ in $\rho_p$ and $\rho_q$ are in ${\mathbb Z}[x]$, pure
of weight $k$
and equal.
\end{cor}
Proof: This follows from lifting both representations simultaneously,
but independently,
and using the Chinese remainder theorem to choose common
characteristic polynomials as we lift using the methods before. 
That the ramified sets for $\rho_p$ and
$\rho_q$  have density $0$ follows from [Kh-Raj]. By construction
we have chosen pure weight $k$ characteristic polynomials 
for all primes outside of the union of the ramified sets
of $\rho_p$ and $\rho_q$ and a finite set.
(We use the larger bounds for the prime $q$ in choosing the sets $R_n$
we need in the Main Theorem.)
The restriction on the weight is necessary if in
$\rhob_p$ inertia at $p$ acts via fundamental characters of level
$2$. See the `Local at $p$ considerations' section of [R3].
\hfill $\square$

\subsection{Some remarks}

Using the methods of this paper, in particular the technique of
proof of the Main Theorem, it is easy to answer in the negative
the following question of [Kh2] (exercise for the interested reader!):
\begin{quest}\label{rationality}
  If $\rho_i:G_{\bf L} \rightarrow GL_m({\bf K})$ 
  is an infinite sequence of (residually absolutely irreducible)
  distinct algebraic representations, all of 
  weight $\leq t$ for some fixed integer $t$, converging 
  to $\rho:G_{\bf L} \rightarrow GL_m({\bf K})$,
  and ${\bf K_i}$ the field of definition of $\rho_i$,
  does $[{\bf K_i}:{\mathbb Q}] \rightarrow \infty$ as $i \rightarrow \infty$?
\end{quest}
Instead of this one should ask the following question to which we
not know the answer:
\begin{quest}
  If $\rho_i:G_{\bf L} \rightarrow GL_m({\bf K})$ 
  is an infinite sequence of (residually absolutely irreducible)
  distinct algebraic representations, 
  such that each $\rho_i$ is {\bf finitely ramified}, and all $\rho_i$ are pure of 
  weight $\leq t$ for some fixed integer $t$, converging 
  to $\rho:G_{\bf L} \rightarrow GL_m({\bf K})$,
  and ${\bf K_i}$ the field of definition of $\rho_i$,
  does $[{\bf K_i}:{\mathbb Q}] \rightarrow \infty$ as $i \rightarrow \infty$?
\end{quest}
The Main Theorem also answers Question 3 of [Kh2] negatively
as using it one gets examples of $\rho$ that are algebraic but infinitely ramified.
Corollary \ref{comp} above also answers negatively Question 1 of [Kh-Raj].

\section{Growth of ramified primes in semisimple
$p$-adic Galois representations}

Let $\rho:G_{\bf K} \rightarrow GL_n({\bf F})$ be a continuous semisimple representation
with ${\bf K}$ a number field and ${\bf F}$ a finite extension of ${\mathbb Q}_p$.
Then it was proven in [Kh-Raj] that the density of the set of primes that
ramify in $\rho$ exists and is 0. For this result the hypothesis that $\rho$ be
semisimple
is crucial. Given any set of primes $T$ 
of ${\bf K}$ that contains all primes above $p$,
using Kummer theory it is easy to construct a non-semisimple 2-dimensional
representation of $G_{\bf K}$ that is ramified at exactly the primes in $T$.

We make a few remarks about [Kh-Raj], and recall some of its results
which suggest Definition~\ref{charpoly}  below.

\begin{prop}
Let $\rho:G_{\bf K} \rightarrow GL_n({\bf F})$ be a continuous semisimple representation
with ${\bf K}$ a number field and ${\bf F}$ a finite extension of ${\mathbb Q}_p$.
Then for all but finitely many primes $q$ of ${\bf K}$, the image of the
decomposition group $D_q$ at $q$ can be conjugated into upper
triangular matrices, and the image of the inertia group $I_q$
at $q$ is unipotent.
\end{prop}
\noindent{Proof:} By Corollary 2 of [Kh-Raj], the image of inertia
group $I_q$ at $q$ is unipotent for almost all primes $q$,
and a fortiori tamely ramified if
$q$ does not lie above $p$ which we now assume. The other assertions
are
contained in the paragraph above Lemma 2 of [Kh-Raj]
which we make more explicit
here. For primes $q$ as in the first sentence, $\rho|_{D_q}$ factors
through
the Galois group of the maximal tamely ramified extension of ${\bf K_q}$,
which is generated by elements $\sigma_q$ and $\tau_q$, such that
$\sigma_q\tau_q\sigma_q^{-1}=\tau^t$ where $t$ is the cardinality of
the residue field at $q$, $\tau_q$ generates tame inertia, and
$\sigma_q$ induces the Frobenius on residue fields. From this relation it follows that
$\sigma_q$ preserves the kernel of $(\tau_q-1)^i$ for any $i$, and
thus as $\rho(\tau_q)$ is unipotent
it follows easily that $\rho(D_q)$ can be conjugated into upper
triangular matrices.
\hfill $\square$

\begin{definition}\label{charpoly}
Let $\rho:G_{\bf K} \rightarrow GL_n({\bf F})$ be a continuous 
semisimple representation
with ${\bf K}$ a number field and ${\bf F}$ a finite extension of ${\mathbb Q}_p$.
Then by the above proposition, for all but finitely many 
primes $q$ of ${\bf K}$, the
image of a
decomposition group $D_q$ at $q$ can be conjugated into upper
triangular matrices, and the image of an inertia group $I_q$
at $q$ is unipotent. For such primes  $q$ we define the characteristic
polynomial of Frobenius $f_q(X)$ at $q$ for the representation $\rho$ to be the
characteristic polynomial of
the image of $\rho$ of any element of $D_q$ which induces the
Frobenius on residue fields.
\end{definition}
\begin{definition}\label{purity}
Let $\rho:G_{\bf K} \rightarrow GL_n({\bf F})$ be a 
continuous semisimple representation
with ${\bf K}$ a number field and ${\bf F}$ a finite 
extension of ${\mathbb Q}_p$. We say that
$\rho$ is density $1$
rational over a number field ${\bf L}$ 
if for a density one set of primes $\{r\}$ at which $\rho$ is unramified,
the characteristic polynomial attached to the conjugacy class
of ${\rm Frob}_r$, the Frobenius 
substitution at $r$, under $\rho$ has coefficients in ${\bf L}[X]$. If further for
for a density one set of primes $\{r\}$ at which $\rho$ is unramified,
the characteristic polynomial attached to the conjugacy class
of ${\rm Frob}_r$ under $\rho$ has roots that are Weil numbers of weight $k$ we
say that $\rho$ is density $1$ pure of weight $k$.
\end{definition}
\begin{definition}\label{compatible}
Let $\rho:G_{\bf K} \rightarrow GL_n({\bf F})$, 
and  $\rho:G_{\bf K} \rightarrow GL_n({\bf F'})$ be
continuous semisimple representations
with ${\bf K}$ a number field and ${\bf F}$ a finite 
extension of ${\mathbb Q}_p$ and ${\bf F'}$ a finite 
extension of ${\mathbb Q}_q$
with $p,q$ primes and fixed embedding of a number field 
${\bf L}$ in ${\bf F}$ and ${\bf F'}$. We
say that $\rho,\rho'$ are density $1$ compatible
if for a density 1 set one set of primes at 
which $\rho,\rho'$ are unramified,
the characteristic polynomial attached to the conjugacy class
of ${\rm Frob}_r$, the Frobenius substitution at $r$, 
in $\rho,\rho'$ are in ${\bf L}[X]$ and equal.
\end{definition}
The following result follows easily from [Kh-Raj] and [S1].
\begin{prop}\label{rub}
  Let $\rho:G_{\bf K} \rightarrow GL_n({\bf F})$ be a 
continuous semisimple representation, and let $G$ be its image
which is a compact $p$-adic Lie group say 
of dimension $N$. Let $C$ be closed subset of 
$G$ which is stable under conjugation
and of $p$-adic analytic dimension $<N$. Then 
the density of the set that consists of primes which are either
ramified in $\rho$, or are unramified
in $\rho$ and such that image of the conjugacy 
class of their Frobenius in $G$ lands in $C$, is 0.
\end{prop}

Serre had asked the first named author 
in an e-mail message in 2000 if the
result of
[Kh-Raj] could be refined to get a stronger 
quantitative control of the growth of primes that ramify in $\rho$.
Serre's question may be motivated by recalling 
the result (Th\'eor\`eme 10 in [S1]) where he shows that
if $\rho$ is finitely ramified, and $C$ is a closed subset of the image of $\rho$ that is stable under conjugation and
of $p$-adic analytic dimension smaller that the $p$-adic 
analytic dimension of the
image of $\rho$,
then \begin{equation}
\pi_C(x)=O\left({x \over {{\rm log}(x)^{1+\epsilon}}}\right)\end{equation}
for some $\epsilon$ bigger than 0  and
with $\pi_C(x)$
the number of primes of norm $\leq x$ whose Frobenius conjugacy class in the image
of $\rho$ lies in $C$. (In fact
under the GRH, Serre proves the stronger estimate $\pi_C(x)=O(x^{1-\epsilon})$ for
some $\epsilon$ bigger than 0.)
In [Kh-Raj] it was shown that
the characteristic polynomials $f_q(X)$ for almost all primes $q$  that
ramify in $\rho$ (in
the sense of Definition \ref{charpoly}) lie in a subvariety of smaller dimension than the character variety of $\rho$.
This together with the results of [S] might lead one to
expect quantitative refinements of the density 0 result of [Kh-Raj]
that better control the order of growth of ramified primes. In fact this is not
the case,
as for example one can prove:

\begin{theorem}\label{growthrate}
  There is a continuous semisimple representation 
  $\rho:G_{\mathbb Q} \rightarrow GL_2({\mathbb Z}_p)$
  (for a prime $p\geq 5$) such that the counting function
  $\pi_{{\rm Ram}(\rho)}(x)$ is not
  $O\left({x \over{{\rm log}(x)^{1+\epsilon}}}\right)$ for any
  $\epsilon>0$, with $\pi_{{\rm Ram}(\rho)}(x)$ the number of
  primes less than $x$ that are ramified in $\rho$.
\end{theorem}

\noindent
We make some preparations for the proof and in particular prove
Lemma~\ref{aux} which is crucial for us.
Consider a surjective
mod $p\geq 5$ Galois representation $\rhobar:G_{\mathbb Q}
\rightarrow GL_2({\mathbb Z}/p{\mathbb Z})$ that
arises from $S_2(\Gamma_0(N))$ for some $(N,p)=1$
with $N$ squarefree. 
Such $\rhobar$ are
plentiful:
for instance take a semistable elliptic curve $E$ over ${\mathbb Q}$ (thus
in particular it does not have CM).
Then for all but finitely many primes $p$, by results of Serre and Wiles, the
corresponding mod $p$ representation will satisfy the
above conditions. Note that as $p>3$ it follows from
a result of Swinnerton-Dyer that any lift $G_{\mathbb Q}
\rightarrow GL_2({\mathbb Z}/p^n{\mathbb Z})$ of $\rhobar$ has image that contains
$SL_2({\mathbb Z}/p^n{\mathbb Z})$.
Let $S$ be the set of ramification of $\rhobar$
(which includes $p$). We have Lemma~\ref{aux}
below that is very close to Lemma 1
of [Kh2]: we give a proof for convenience, as it is crucial to the proof of our
proposition, and because
the last part of it is not covered in [Kh2].

We say that a finite set  $Q$ of nice primes
is {\em auxiliary} if certain maps on $H^1$ and
$H^2$, namely 
\begin{equation}H^1(G_{S \cup Q},{\rm Ad}^0(\rhobar))
\rightarrow \oplus _{v \in S \cup Q}H^1(G_v,{\rm Ad}^0(\rhobar))/{\cal
N}_v
\end{equation}
and 
\begin{equation}
H^2(G_{S \cup Q},{\rm Ad}^0(\rhobar))
\rightarrow \oplus _{v \in S \cup Q} H^2(G_v,{\rm Ad}^0(\rhobar))
\end{equation}
considered in [R3]  are
isomorphisms. We refer to [R3]
for the notation used: recall that ${\cal N}_v$ for $v \in  Q$
is the mod $p$ cotangent space of a smooth quotient of the local
deformation ring at $v$ which parametrises lifts whose mod $p^m$
reductions are $\rho_m$-nice for all $m$.
Henceforth we call such representations
{\em special}.
Here for $v \neq p$ and
$v \in S$, ${\cal N}_v$ is described as the image
under the inflation map of
$H^1(G_v/I_v,\Ad^0(\rhobar)/{\Ad^0(\rhobar)}^{I_v})$.
For $v=p$ we define it to be either 
$H^1_{fl}(G_p,\Ad^0(\rhobar))$ or
$H^1_{Se}(G_p,\Ad^0(\rhobar))$
according to whether
$\rhobar|_{I_p}$ is finite or not, using the notation of Section 4.1
of de Shalit's article in [FLT].
For $v \in Q$, ${\cal N}_v$ is described as the
subspace of $H^1(G_v,\Ad^0(\rhobar))$ generated by the class of the
cocycle that (in a suitable choice of basis of $V_{\rhobar}$, the
2-dimensional $k$-vector space that affords $\rhobar$, and viewing
$\Ad^0(\rhobar)$ as a subspace of ${\rm End}(V_{\rhobar})$)
sends $\sigma \rightarrow 0$ and $\tau$ to
\begin{equation}\left(\matrix{0&1\cr
		0&0}\right)\end{equation}
where $\sigma$ and $\tau$ generate the tame quotient of
$G_q$, and satisfy the relation $\sigma\tau\sigma^{-
1}=\tau^q$.

The isomorphisms above result in the fact that there is a
lift $\rho_{S \cup Q}^{Q-new}:G_{\mathbb Q} \rightarrow GL_2(W(k))$ 
of $\rhobar$
which is furthermore the unique lift of $\rhobar$ to a representation
to $GL_2({\cal O})$ (with ${\cal O}$ ring of integers of
any finite extension of ${\mathbb Q}_p$)
that has the properties of being semistable
of weight 2, unramified outside $S \cup Q$, minimally ramified at
primes in $S$, with determinant $\varepsilon$, and {\it special} at primes in $Q$.

\begin{lemma}\label{aux}
Let $\rho_n:G_{\mathbb Q} \rightarrow
GL_2({\mathbb Z}/p^n{\mathbb Z})$ be any lift of $\rhobar$ that is minimally
ramified at the primes in $S$ and such that
all other primes that ramifiy in $\rho_n$ are special and
nice. Let $Q_n'$ be any finite set primes
containing the ramification of $\rho_n$
such that $Q_n' \backslash S$ contains only $\rho_n$-nice primes.
Then there exists a finite set of primes $Q_n$ that contains $Q_n'$,
such that $\rho_n|_{D_q}$ is special for $q \in Q_n \backslash S$, $Q_n \backslash
S$ contains only nice primes
and $Q_n \backslash S$ is auxiliary. Further the representation
$\rho_{Q_n}^{Q_n \backslash S-new}:G_{\mathbb Q} \rightarrow GL_2(W(k))$ is ramified at all primes in
$Q_n$, and mod $p^n$ is isomorphic to $\rho_n$.
\end{lemma}

\noindent{Proof of Lemma:} 
We use [R3] and Fact~\ref{disjointness} (that latter
being a certain mutual disjointness result for field extensions cut out by
$\rho_n$ and extensions cut out by  elements of
$H^1(G_{\mathbb Q},{\rm Ad}^0(\rhobar))$
and $H^1(G_{\mathbb Q},\adst)$ 
and here we use $p>3$)
to construct an auxiliary set of primes $V_n$
such that $\rho_n|_{D_q}$ is special for $q \in V_n$. Then
as $Q_n' \backslash S$ contains only
nice primes, it follows from Proposition 1.6 of [W]
that the kernel and cokernel of the map
\begin{equation}H^1(G_{S \cup V_n \cup Q_n'},\ad)
\rightarrow \oplus _{v \in S \cup V_n \cup Q_n'}
H^1(G_v,\ad)/{\cal N}_v\end{equation}
have the same cardinality, as the domain and range have the same
cardinality. Then using Proposition 10 of [R3],
or Lemma 1.2 of [T], and Lemma 8 of [Kh-Ram],
we can augment the set $S \cup V_n \cup Q_n'$ to
get a set $Q_n$ as in the statement of the lemma. The last line follows from
the fact that by Theorem 1 of [Kh1]
$\rho_{Q_n}^{Q_n-new}$ arises from a newform (note that the $\rhobar$ we have fixed is modular!), and thus is forced
to be ramified at primes in $Q_n$
for purity reasons. Further,
$\rho_{Q_n}^{Q_n-new}$ mod $p^n$ is isomorphic to $\rho_n$
as  there is a unique mod $p^n$ lift of
$\rho$ that is minimally ramified at primes
in $S$ and special for all primes in
$Q_n \backslash S$ with a certain fixed determinent character.
As $\rho_n$ and $\rho_{Q_n}^{Q_n\backslash S-new}$ mod $p^n$ are two such
lifts they are forced to be isomorphic.
\hfill $\square$

\vskip1em\noindent
Proof of Theorem \ref{growthrate}:
We construct the $\rho$ of the theorem as the inverse limit of a
compatible family of mod $p^n$ representations $\rho_n$ for an
infinite number of positive integers $n$.
When lifting $\rhobar$ to a mod $p^{m_1}$ representation (for some large $m_1$)
we
choose an integer $f_1$ large enough
so that there are at least
${{f_1} \over {{\rm log}({f_1})^{2}}}$ nice primes up to $f_1$.
This can be done as
the set of nice primes has
positive density by the Cebotarev density theorem. Apply the lemma
with $n=1$, taking $Q_1'$ to be set that contains $S$
and all the nice primes up to $f_1$. Thus we get
an auxiliary set $Q_1$ that contains $Q_1'$.
Consider $\rho_{Q_1}^{Q_1\backslash S-new}$, and an integer
$m_1$ such that $\rho_{Q_1}^{Q_1\backslash S-new}$
mod $p^{m_1}$ is ramified at all primes in $Q_1$: such an integer $m_1$ exists
because of the last line of the lemma.
When lifting $\rho_{Q_1}^{Q_1\backslash S-new}$ mod $p^{m_1}$ to a mod $p^{m_2}$ representation (for some large $m_2$)
we choose an integer $f_2\gg f_1$ large enough
so that up to $f_2$ there are at least
${{f_2} \over {{\rm log}(f_2)^{3/2}}}$
$\rho_{Q_1}^{Q_1\backslash S-new}$ nice primes:
again by the Cebotarev density theorem (and the largeness of
the image of $\rho_{Q_1}^{Q_1\backslash S-new}$ mod $p^{m_1}$)
all large enough $f_2$ will satisfy this property.
Choose $t_2'$
to be the set that contains $S$ and all the nice primes that
are special for $\rho_{Q_1}^{Q_1\backslash S-new}$
mod $p^{m_1}$ up to $f_2$, together with all the primes in $Q_1$.
Applying the lemma we get
an auxiliary set $Q_2$ that contains $Q_2'$ such that all primes in $Q_2 \backslash S$ are
special	 for $\rho_{Q_1}^{Q_1\backslash S-new}$ mod $p^{m_1}$. By the last line
of
the lemma,
because of this $\rho_{Q_2}^{Q_2\backslash S-new}$ mod $p^{m_1}$ is isomorphic to $\rho_{Q_1}^{Q_1\backslash S-new}$ mod $p^{m_1}$.
Further there is a $m_2\gg 0$ such that $\rho_{Q_2}^{Q_2\backslash S-new}$ mod $p^{m_2}$
is ramified at all primes in $Q_2$. Now the inductive procedure is clear.
At the $n$th stage, we will be dealing with $\rho_{Q_n}^{Q_n\backslash S-new}$,
and an integer
$m_n$ such that $\rho_{Q_n}^{Q_n\backslash S-new}$ mod $p^{m_n}$ is ramified at
all primes in $Q_n$.
When lifting $\rho_{Q_n}^{Q_n\backslash S-new}$ mod $p^{m_n}$ to a mod $p^{m_{n+1}}$ representation (for some large $m_{n+1}$)
we choose an integer $f_{n+1}\gg f_n$ large enough
so that up to $f_{n+1}$ 
there are at least ${ {f_{n+1}} \over
{ {\rm log}(f_{n+1})^{1+{1 \over {2^n}}}}}$ 
nice primes that are special for $\rho_{Q_n}^{Q_n\backslash S-new}$
mod $p^{m_n}$.
By the Cebotarev density theorem all large enough $f_{n+1}$
will satisfy this property. Choose $Q_{n+1}'$
to be a set that contains $S$, contains $Q_n$ and contains
all the nice primes that are special for $\rho_{Q_n}^{Q_n\backslash S-new}$ mod
$p^{m_n}$ up to $f_{n+1}$.
Applying the lemma we get
an auxiliary set $Q_{n+1}$ that contains $Q_{n+1}'$ such that all primes in $Q_{n+1} \backslash S$ are
special	 for $\rho_{Q_{n}}^{Q_n\backslash S-new}$ mod $p^{m_n}$. Thus by the last line of lemma,
$\rho_{Q_{n+1}}^{Q_{n+1}\backslash S-new}$ mod $p^{m_{n}}$ is isomorphic to $\rho_{Q_{n}}^{Q_{n}\backslash S-new}$ mod $p^{m_n}$.
Further there is a $m_{n+1}\gg 0$ such that
$\rho_{Q_{n+1}}^{Q_{n+1}\backslash S-new}$ mod $p^{m_{n+1}}$
is ramified at all primes in $Q_{n+1}$. In this process we also take care (as we may easily do!) to choose the sequences $m_n$ and
$f_n$ so that they tend to infinity with $n$ (in fact if the $f_n$'s tend to infinity,
so do the $m_n$'s). We define $\rho$ to be the inverse limit
of the compatible system of
representations $\rho_{Q_{n}}^{Q_{n}\backslash S-new}$ mod $p^{m_n}$.
\hfill
$\square$

\vspace{3mm}

\noindent{Remarks:}
1) It is interesting to note that although the
$\rho$ we construct is infinitely ramified, it is important for us
to construct it as a limit of {\it geometric} representations $\rho_i$
which are in particular finitely ramified, as the geometricity of the
$\rho_i$'s is vital in ensuring that the limit $\rho$ is ramified at
very many primes!
\newline\noindent
2) The proposition is not the best possible result and merely
   illustrates the fact that growth of ramified primes can be rapid. By the same   
methods as used in the proof above we can construct a semisimple  $\rho:G_{\mathbb Q} \rightarrow GL_2({\mathbb Z}_p)$
 such that for a given $n$ the counting function
  $\pi_{{\rm Ram}(\rho)}(x)$ is not
  $O(x/{\rm log}(x)({{\rm log}^{(n)}(x))}^{\epsilon})$ for any
  $\epsilon>0$ where ${\rm log}^{(n)}$
means ${\rm log}$ composed with itself $n$ times.

\section*{References}

\noindent [B] B\"{o}ckle, G.,
{\it A local-to-global principle for deformations of Galois representations},
J. Reine Angew. Math. 509 (1999), 199--236.

\vspace{3mm}

\noindent [DT] Diamond, F., Taylor, R., 
{\it Lifting modular mod $l$ representations},
Duke Math. J. 74 (1994), no. 2, 253--269.

\vspace{3mm}

\noindent [FLT] {\it Modular forms and Fermat's last theorem}, 
edited by Gary Cornell, Joseph H. Silverman and Glenn Stevens.
Springer-Verlag, New York, 1997.

\vspace{3mm}

\noindent [Kh1] Khare, C., {\it On isomorphisms between deformation rings and
Hecke rings}, with an appendix by G. B\"ockle, 
Invent. Math. 154 (2003), vol. 1, 199 - 222.

\vspace{3mm}

\noindent [Kh2] Khare, C., {\it Limits of residually irreducible
$p$-adic Galois representations}, Proc. Amer. Math. Soc. 131 (2003), 1999-2006.

\vspace{3mm}

\noindent [Kh3] Khare, C., {\it $F$-split Galois representations are potentially abelian}, Proc. Amer. Math. Soc. 
131 (2003), 3021-3023.

\vspace{3mm}

\noindent [Kh-Raj] Khare, C., Rajan, C.~S.~, {\it The density
of ramified primes in semisimple $p$-adic Galois representations},
International Mathematics Research Notices  no. 12 (2001), 601--607.

\vspace{3mm}

\noindent [Kh-Ram] Khare, C., Ramakrishna, R., {\it
Finiteness of Selmer groups and deformation rings}, Invent. Math. 154 (2003), vol. 1, 179--198.

\vspace{3mm}
\noindent [M] Mazur, B., {\it An introduction to the deformation theory of Galois representations}, in `Modular forms and Fermat's last theorem',
edited by Gary Cornell, Joseph H. Silverman and Glenn Stevens.
Springer-Verlag, New York, 1997.

\vspace{3mm}

\noindent [NSW] Neukirch, J., Schmidt, A., Wingberg, K., {\it Cohomology of number fields}, Grundlehren der Mathematischen
   Wissenschaften 323, Springer-Verlag, 2000.

\vspace{3mm}

\noindent [Ri] Ribet, K., {\it Report on mod $\ell$
representations of ${\rm Gal}(\overline{\mathbb Q}/{\mathbb Q})$}, in Motives, Proc.
Sympos. Pure Math. 55, Part 2 (1994), 639--676.

\vspace{3mm}
\noindent [R1] Ramakrishna, R., {\it Lifting Galois representations},
 Inventiones. Math. 138 (1999), no. 3, 537--562. 

\vspace{3mm}

\noindent [R2] Ramakrishna, R., {\it Infinitely ramified representations},
Annals of Mathematics 151 (2000), 793--815.

\vspace{3mm}

\noindent [R3] Ramakrishna, R., {\it Deforming Galois representations and the
conjectures of Serre and Fontaine-Mazur}, Annals of Mathematics 156 (2002), 115--154.

\vspace{3mm}

\noindent [S1] Serre, J-P., {\it  Quelques applications du theoreme de
densite de Chebotarev}, Collected Works, Vol. 3, 563--641.

\vspace{3mm}
                                                                                
\noindent [T1] Taylor, R., {\it Remarks on a conjecture of 
Fontaine and Mazur}, J. Inst. Math. Jussieu  1  (2002),  no. 1, 125--143.

\vspace{3mm}

\noindent [T2] Taylor, R., {\it On icosahedral Artin
representations II},  Amer. J. Math.  125  (2003),  no. 3, 549--566.

\vspace{3mm}

\noindent[TW] Taylor, R., Wiles, A., 
{\it Ring-theoretic properties of certain Hecke
algebras}, Annals of Mathematics (2) 141 (1995), 553--572. 

\vspace{3mm}

\noindent [W] Wiles, A., {\it Modular elliptic curves and 
Fermat's last theorem}, Annals of Mathematics 141 (1995), 443--551. 

\vspace{3mm}

\noindent {\it Addresses:}

\noindent CK: 155 S 1400 E, Dept of Math, Univ of Utah, Salt Lake City, UT 84112, USA.
e-mail: {\tt shekhar@math.utah.edu}

\vspace{3mm}

\noindent ML:  Department of Mathematics, Indiana University, Bloomington, IN 47405, 
USA. e-mail: {\tt larsen@iu-math.math.indiana.edu }

\vspace{3mm}

\noindent RR: Department of Mathematics, Cornell University, Malott Hall,
Ithaca, NY 14853, USA. e-mail: {\tt ravi@math.cornell.edu}

\end{document}